\documentclass[11pt,reqno]{amsart}
\usepackage{amsmath,amsthm,amssymb,latexsym, eucal}

\newtheorem{thm}[equation]{Theorem}
\newtheorem{cor}[equation]{Corollary}
\newtheorem{lem}[equation]{Lemma}

\newtheorem{exams}[equation]{Examples}
\newtheorem{claim}[equation]{Claim}
\newtheorem{rem}[equation]{Remark}

\newtheorem{subsec}[equation]{}
\theoremstyle{definition}
\newtheorem{defn}[equation]{Definition}

\numberwithin{equation}{section}

\newcommand{\ot}{\otimes}

\newcommand{\hf}{\frac{1}{2}}
\newcommand{\pf} {\medskip \noindent \textbf{Proof.\  \  }}
\newcommand{\supp}{\hbox{\rm supp\,}}

\newcommand{\al}{\alpha}

\newcommand{\De}{\Delta}

\newcommand{\m}{\medskip}
 
\newcommand{\tr}{\mathfrak{tr}}

\newcommand{\vvp}{\varphi}

\newcommand{\F}{\mathbb F}

\newcommand{\fa}{{\mathfrak a}} 
\newcommand{\g}{{\mathfrak g}} 
\newcommand{\h}{{\mathfrak h}}
\newcommand{\fs}{{\mathfrak s}} 

\newcommand{\scl}{{\mathcal L}}

\newcommand{\Id}{\mathop{\rm{ id}}}

\newcommand{\inder}{\mathop{\rm IDer}}

\begin{document}

\title[]{\large Lie $\boldsymbol{G}$-Tori of  Symplectic Type}

\author[]{Georgia Benkart$^{1,*}$}
\thanks{$^{1}$ Supported in part by National Science Foundation Grant
\#{}DMS--0245082.} 
\author[]{Yoji Yoshii$^{2}$}
\thanks{$^{2}$ Support  from an NSERC postdoctoral fellowship during time at the
University of Wisconsin-Madison  is gratefully acknowledged.
\newline{$^{*}$Corresponding author}
} 
\date{\today} 
\begin{abstract}  We classify centerless Lie
$G$-tori of type $\hbox{\rm C}_r$ including the most difficult case $r=2$ by applying techniques
due to Seligman. 
In particular, we show that the coordinate algebra of a Lie $G$-torus of type $\hbox{\rm C}_2$ is either an associative $G$-torus with involution or a Clifford $G$-torus.      Our results generalize the 
classification of the core of the extended affine Lie algebras of type $\hbox{\rm C}_r$ by Allison and Gao.   \end{abstract}
\maketitle

\begin{center}{\small  {\it Dedicated to Professor George Seligman
with admiration}}\end{center}

\section{Introduction}

The extended affine Lie algebras of \cite{AABGP} are natural generalizations of the affine and 
toroidal Lie algebras,  which have played such a pivotal  role in 
diverse areas of mathematics and physics.   Their root 
systems (the so-called extended affine root systems) are essential
in the work of Saito (\cite{S1}, \cite{S2}) and Slodowy 
\cite{Sl} on singularities.    An extended affine Lie algebra $\mathcal E$  possesses a nondegenerate 
invariant symmetric
bilinear form and a finite-dimensional self-centralizing ad-diagonalizable subalgebra $\mathcal H$.
The  root system $R$ of $\mathcal E$  relative to $\mathcal H$ has a decomposition  $R = R^0 \cup R^\times$  into  
isotropic roots  $R^0$  and  nonisotropic roots $R^\times$.     The
subalgebra $\mathcal E_c$  generated by the root spaces corresponding to the nonisotropic roots 
is an ideal of $\mathcal E$,  called the {\em core} of $\mathcal E$.   The algebra $\mathcal E$
is said to be {\em tame} if the centralizer of $\mathcal E_c$ in $\mathcal E$ is just  the center ${\mathcal Z}(\mathcal E_c)$ of $\mathcal E_c$. 

The core features prominently in the classification of the tame extended affine
Lie algebras (see \cite{AABGP}, \cite{BGK}, \cite{BGKN}, \cite{AG},  \cite{AY}, \cite{ABG3}).   It
is graded by a finite (possibly nonreduced) irreducible root system $\Delta$   and   
is root-graded in the sense of \cite{BM} and \cite{ABG2}.    It also has a grading by the
free abelian group $\Lambda$ generated by $R^0$.   Moreover,  $\mathcal E_c/
{\mathcal Z}(\mathcal E_c)$  is what is now referred to as a centerless  Lie  torus (as in \cite{N1} and 
\cite{N2}), and every centerless Lie torus is the centerless  core of an extended affine Lie algebra  (see \cite{Y6}).
The family of centerless Lie tori consists of mostly
infinite-dimensional Lie algebras, but these Lie algebras are well-behaved counterparts of
the finite-dimensional split simple Lie algebras.    
Roughly speaking, the axioms
that characterize tame extended affine Lie algebras and that get
used heavily in the determination of their cores are replaced in the notion of a Lie torus
by gradings by an irreducible root system $\Delta$  and by a free abelian group $\Lambda$  of finite rank,  by a certain property that is sometimes called the division
property, and by the assumption that the homogeneous pieces have dimension at most 1.  
When the rank of $\Lambda$ is $\geq 1$, the Lie algebra is infinite-dimensional, but when it
is 0, the algebra is just a finite-dimensional split simple Lie algebra. 
 
In this paper,  we study Lie $G$-tori, where the free abelian group $\Lambda$ in
the definition of a Lie torus  is replaced
by an arbitrary abelian group $G $.  
Lie $G$-tori were first introduced by the second author in \cite{Y5} and  \cite{Y6} as a special class of root-graded Lie algebras.   
The classification of Lie $G$-tori of type $\text A_r$ can be easily
derived  from  results  in \cite{BGK},  \cite{BGKN}, \cite{Y1}, \cite{Y2}, \cite{Y4} or \cite{AY}.
They are coordinatized by $G$-tori (in the sense of Definition \ref{def:Gtorus} below), 
which are associative when $r \geq 3$,
alternative  when $r=2$, and  Jordan  
when $r=1$.   The Lie $G$-tori of type $\text B_r$ ($r \geq 3$) were determined
in \cite{Y5}, where it was shown that their coordinate algebras are a special type of Jordan $G$-torus,
called a Clifford $G$-torus.      Our aim is to classify centerless  Lie $G$-tori of type $\text B_2=\hbox{\rm C}_2$  together with those of type  $\hbox{\rm C}_r$ for higher rank $r \geq 3$.
Lie $G$-tori of type $\text F_4$, $\text G_2$, and $\text{BC}_r$  have not yet been classified,  although the centerless cores of  extended affine Lie algebras 
(the centerless Lie tori) of types $\text F_4$,  $\text G_2$, and $\text{BC}_r$   have been 
determined  in \cite{AG},  \cite{AY},  \cite{ABG3}, \cite{AFY}, and \cite{F}.   
 \bigskip
 \section{ Preparation}
 \m
 {\em Throughout we will assume that all algebras are over a field $\F$ of characteristic zero.} 
 \m

Let $\Delta$ be a finite irreducible root system (not necessarily reduced) as in
\cite[Chap.~VI, \S 1.1]{B}.   For  each root $\mu \in \Delta$, let  $\mu^\vee$ be the corresponding coroot so that 
$\langle \nu,\mu^\vee\rangle = 2(\nu,\mu)/(\mu,\mu)$ is the Cartan integer
for all $\nu \in \Delta$.      Set $\Delta_{\text {ind}} = \{ \mu \in \Delta \mid \hf \mu \not \in \Delta \}$.   Let $G = (G,+,0)$ be an additive abelian group.    For any subset $S$ of $G$, 
we denote the subgroup of $G$ that $S$ generates  by $\langle S \rangle$.
\medskip

\begin{defn}\label{defn:LGT}  A Lie algebra $\scl$  is a {\em Lie 
$G$-torus of type $\Delta$}  if 
\begin{itemize}
\item[{\rm (1)}]  $\scl$ has a decomposition into subspaces  $\scl = \displaystyle{\bigoplus_{\mu \in \Delta \cup \{0\}, \, g \in G}}   \scl^g_\mu$ such that  $[\scl^g_\mu, \scl^h_\nu]  \subseteq \scl^{g+h}_{\mu+\nu}$;
\item[{\rm (2)}]   For every $g \in G$, \ $\scl_0^g = \displaystyle{\sum_{\mu \in \Delta,
h \in G}} [\scl_\mu^h, \scl_{-\mu}^{g-h}]$; 
\item[{\rm (3)}]  

\begin{itemize}
\item[{\rm (a)}]  For each nonzero $x \in \scl_\mu^g$ ($\mu \in \Delta, g \in G)$,  there exists a $y \in 
\scl_{-\mu}^{-g}$  so that 
$t : = [x, y] \in \scl_0^0$ satisfies 
$[t, z] = \langle \nu,\mu^\vee\rangle z$ for all  $z \in \scl_\nu^h$, 
$(\nu \in \Delta \cup \{0\}, h \in G)$. 
\item[{\rm (b)}]  $\dim \scl_\mu^g \leq 1$ and $\dim \scl_\mu^0 = 1$ if 
$\mu \in \Delta_{\small{\textrm{ind}}}$;
  \end{itemize}
\item[{\rm (4)}]  $G = \langle \hbox{\rm supp}\,\scl \rangle$ where 
$$ \hbox{\rm supp}\,\scl := \{g \in G \mid \scl_\mu^g \neq 0 \ \hbox{\rm  for some} \ \mu \in \Delta
\cup \{0\}\}.$$
\end{itemize}
\end{defn}

\begin{subsec}\label{rems:2.2} {\rm \textbf{Remarks on Definition \ref{defn:LGT}}} \end{subsec}

Condition (4) is simply a convenience. If it fails to hold, we may replace
$G$ by the subgroup generated by $ \hbox{\rm supp}\,\scl$. \medskip

It follows from (1) that $\scl$ is graded by the group $G$.  Thus, if
$\scl^g := \bigoplus_{\mu \in \Delta \cup \{0\}} \scl_\mu^g$, then
$\scl = \bigoplus_{g \in G} \scl^g$ and $[\scl^g,\scl^h] \subseteq \scl^{g+h}$.  \medskip

The Lie algebra $\scl$ also admits a grading by the root lattice $Q(\Delta)$:  if 
$\scl_\lambda := \bigoplus_{g \in G}  \scl_\lambda^g $ for 
$\lambda \in Q(\Delta)$, where $\scl_\lambda^g = 0$ if $\lambda \not \in \Delta \cup \{0\}$, then $\scl = \bigoplus_{\lambda \in Q(\Delta)} \scl_\lambda$
and  $[\scl_\lambda, \scl_\mu] \subseteq \scl_{\lambda+\mu}$. \medskip

 {F}rom (3) we see for $\mu \in \Delta_{\textrm {ind}}$ that there
 exist  elements $e_\mu \in \scl_\mu^0$, $f_\mu \in \scl_{-\mu}^0$, 
 and $t_\mu: = [e_\mu,f_\mu]$ so that $[t_\mu, z] = \langle \nu,\mu^\vee \rangle z$
 for all $z \in \scl_\nu^h$, 
$(\nu \in \Delta, h \in G)$.  Thus,  the elements
$e_\mu,f_\mu,t_\mu$ determine a canonical basis for
a copy of the Lie algebra $\mathfrak{sl}_2$.     In addition,  the products $[t_\lambda, t_\mu] = 0$  for
$\lambda,\mu \in  \Delta_{\textrm {ind}}$.  
 It follows that  the subalgebra 
$\g$ of $\scl$ generated by the subspaces $\scl_\mu^0$ for
$\mu \in \Delta_{\textrm {ind}}$  is a split simple Lie algebra with
split Cartan subalgebra $\mathfrak h: = \sum_{\mu \in \Delta_{\textrm {ind}}}
[\scl_\mu^0, \scl_{-\mu}^0]$.       We may identify the coroot $\mu^\vee$ with 
an element of $\mathfrak h$.    Then $t_\mu$ is equivalent to $\mu^\vee$ modulo
the center $\mathcal Z(\scl)$ of $\scl$.  
Condition (3a) (or equivalently,  existence of  canonical $\mathfrak{sl}_2$-basis elements $e_\mu
\in \scl^0_\mu, f_\mu \in \scl^0_{-\mu}$, and $t_\mu = [e_\mu,f_\mu]$
such that $t_\mu \equiv \mu^\vee$ mod $\mathcal Z(\scl)$  for each $\mu \in \Delta_{\textrm {ind}}$)
is often called
the {\em division property}, and $\scl$ is said to be {\em division graded}.   \medskip

It follows then that  a  Lie $G$-torus  $\scl$ is a Lie
algebra graded by the root system $\Delta$  in the following sense: 

\begin{defn}\label{rootgrade} A Lie algebra $\scl$  is said to be  {\em graded by the root system  
$\Delta$} (where $\Delta$ is finite and irreducible)   or  to be {\it $\Delta$-graded}  if
\begin{itemize}
 \item[{\rm ($\Delta$1)}]  $\scl$ contains as a subalgebra  a
finite-dimensional split ``simple'' Lie algebra $\g$, called the 
{\it grading subalgebra},  with root system 
$\Delta_\g$ relative to a split Cartan subalgebra $\mathfrak h$; 

 \item[{\rm ($\Delta$2)}] $\scl$ has a decomposition into subspaces $ \scl= \bigoplus_{\mu \in \Delta \cup \{0\}}
 \scl_\al$,  where $\scl_\mu = \{v \in \scl \mid [t,v] = \mu(t)v$ for all $t \in \mathfrak h\}$.   
 \item[{\rm ($\Delta$3)}]  $\scl_0 = \sum_{\mu \in \Delta}\, [\scl_\mu, \scl_{-\mu}]$;
 \item[\rm ($\Delta$4)] either $\Delta$ is reduced and equals the root system
$\Delta_\g$ of $(\g,\mathfrak h)$  or $\Delta={\rm BC}_r$ and $\Delta_\g$ is of type  
B$_r$,  C$_r$,  or D$_r$.
\end{itemize} \end{defn}

The word simple is in quotes above, because in all instances except two, 
$\g$ is a simple Lie algebra.  The sole exceptions are when $\Delta$ 
is of type BC$_2$,   $\Delta_\g$ is of type D$_2$ = A$_1 \times$ 
A$_1$,  and $\g$ is the direct sum of two copies of 
$\mathfrak{sl}_2$;  or when $\Delta$ is of type BC$_1$, $\Delta_\g$ 
is of type D$_1$, and $\g = \mathfrak h$ is one-dimensional.  Neither
of these exceptions will play a role in this work.   \m

 The definition above is due to Berman-Moody \cite{BM} for the 
case $\De=\De_\g$.  The extension to the nonreduced root systems BC$_r$ was 
developed by  Allison-Benkart-Gao in  \cite{ABG2} for $r \geq 2$  and  by  Benkart-Smirnov
in \cite{BS} for $r =1$.       \m 

A Lie $G$-torus $\scl$ of type $\hbox{\rm BC}_r$  has grading subalgebra 
$\g$ generated by the root spaces  $\scl^g_\mu$  for 
$\mu \in \Delta_{\small{\textrm{ind}}}$, and so $\g$ will be
of type $\hbox{\rm B}_r$,  since 
those
root spaces have dimension one by 
 (3b) of Definition \ref{defn:LGT}.      \m

The original definition of a Lie $G$-torus in \cite{Y5} required the Lie algebra $\mathcal L$
to be $\Delta$-graded.    As mentioned above, this holds automatically.    
\medskip

Sometimes in what follows, we stipulate that a Lie algebra
is {\em $(\Delta,G)$-graded}.    By that we mean it is a $\Delta$-graded Lie algebra $\scl$
which is also $G$-graded, $\scl = \bigoplus_{g \in G}\scl^g$,  such that 
the grading subalgebra $\g$ of $\scl$  is contained in $\scl^0$ and 
the support $\{g \in G \mid \scl^g \neq 0\}$ generates $G$.   It follows that every 
$(\Delta,G)$-graded Lie algebra has a decomposition $\scl = \displaystyle{\bigoplus_{\mu \in \Delta \cup \{0\}, \, g \in G}}   \scl^g_\mu$,  where $\scl_\mu^g = \scl_\mu \cap \scl^g$,
and that condition (4) of Definition \ref{defn:LGT} holds.  
 
\m

In \cite[Defn.~3.12]{BN}, a Lie $G$-torus is defined to be a $(\Delta,G)$-graded
Lie algebra $\scl$ satisfying (3) of Definition \ref{defn:LGT}.  The definition
in \cite{BN} permits  the grading subalgebra to be type $\hbox{\rm C}_r$ ($r \geq 1$) or $\hbox{\rm D}_r$ ($r \geq 3$) when $\Delta$ is of type $\hbox{\rm BC}_r$.    \m

\begin{subsec} {\rm \textbf{Lie  algebras graded by C$_r$, $r \geq 2$.}} \end{subsec}

We specialize now to Lie algebras graded by the root systems C$_r$, $r \geq 2$,
since ultimately we intend to classify the centerless Lie $G$-tori of type C$_r$.
Towards this purpose, let 
$V$ be a $2r$-dimensional vector space over  
$\F$  with a nondegenerate skew-symmetric bilinear form.
 Let  $\{v_1,\ldots,v_{2r}\}$ be a basis of $V$,
 and let  $\g$ denote the symplectic Lie algebra $\mathfrak{sp}(V)$ of skew
endomorphisms of $V$.    Using the basis above, we identify
$\g$ with $\mathfrak {sp}_{2r}(\F)$, the Lie algebra of $2r \times 2r$ matrices $x$ over $\F$ that  satisfy $x^{t}M = -Mx$, where $M$ is the
matrix whose $(i,j)$-entry is 
$\text{sign} (i-j)\delta_{i+j, 2r+1}$ for $1 \leq i,j \leq 2r$.  We also identify a Cartan subalgebra $\h$
of $\g$ with the set of diagonal matrices in $\g$.
Thus, the elements  $\{E_{1,1}-E_{2r, 2r},\dots, E_{r,r}-E_{r+1, r+1}\}$
determine a basis for $\h$, where the $E_{i,j}$ are the standard matrix units.  
Let $\{\varepsilon_1,\ldots,\varepsilon_r\}$
be the dual basis in $\h^*$.
Then $\g$ has a decomposition into one-dimensional root spaces relative to $\h$, and 
a basis for these root spaces may be taken as follows:
\begin{itemize}
\item[{($\g$1)}]  $E_{i,j}-E_{2r+1-j, 2r+1-i}$ \  for  \  $\varepsilon_i-\varepsilon_j$,
\item[{($\g$2)}] $E_{i, 2r+1-j} + E_{j, 2r+1-i}$ \  for \   $\varepsilon_i+\varepsilon_j$,  
\item[{($\g$3)}]
$E_{2r+1-j,i} + E_{2r+1-i,j}$ \ for \ $-\varepsilon_i-\varepsilon_j$
\end{itemize}
where  $1 \leq i,j \leq r$.   The corresponding root system $\Delta$ decomposes  into the set 
$\Delta_{sh}:=\{\pm (\varepsilon_i\pm\varepsilon_j)\mid 1\leq i\neq j\leq r\}$ of  short roots 
and the set $\Delta_{lg}=\{\pm 2\varepsilon_i\mid 1\leq i\leq r\}$ of long roots.
\m
 
Let $\fs$ denote the set of  $2r \times 2r$ matrices $s$  of trace zero  satisfying
 $s^{t} M = M s$.   Then $\fs$ is a $\g$-module under the action
 $x.s = [x,s] = xs-sx$  $(x \in \g, s \in \fs)$,   and $\fs$ has a decomposition into one-dimensional 
 weight spaces relative to $\h$.  A basis for these weight spaces may be chosen as follows:
\begin{itemize}
\item[{($\fs$1)}]  $E_{i,j}+E_{2r+1-j, 2r+1-i}$ \ for  
$\begin{cases}\varepsilon_i+\varepsilon_j & \hbox{ \rm if}  \  \  i<j \\
 -\varepsilon_i-\varepsilon_j& \hbox{ \rm if}  \ \   i>j \end{cases}$
\item[{($\fs$2)}] $E_{i, 2r+1-j} - E_{j, 2r+1-i}$ \ for \ $\varepsilon_i-\varepsilon_j$,  
\item[{($\fs$3)}]  $E_{2r+1-m,i} - E_{2r+1-i,j}$ \  for \ $-\varepsilon_i+\varepsilon_j$,
\end{itemize} 
where $1 \leq i\neq j \leq r$.  \m

A  $\hbox{\rm C}_r$-graded Lie algebra $\scl$ decomposes
into copies of $\g$, $\fs$, and the trivial one-dimensional $\g$-module relative
to the adjoint action of the grading subalgebra $\g$.  By collecting isomorphic
summands, we may assume there are $\F$-vector spaces $A,B,D$ so that 
$$\scl=(\g\otimes A)\oplus (\fs \otimes B)\oplus D,$$ 
where  $D$ is the sum of the trivial modules.    By \cite{ABG1}, there is a symmetric product $\circ$ and a skew-symmetric product $[\cdot,\cdot]$
on $\fa = A \oplus B$ so that $\fa$ with the multiplication 
\begin{equation}\label{eq:comult} \alpha\alpha' = \hf(\alpha \circ \alpha') + \hf [\alpha,\alpha']
\end{equation}
for $\alpha,\alpha' \in \fa$ is the {\em coordinate algebra} of $\scl$.
The space $D$ is a Lie subalgebra
of $L$, which acts as derivations on $\fa$.     When $\scl$ is centerless,
then $D$ is spanned by the inner derivations $D_{\alpha,\alpha'}$ for $\alpha,\alpha' \in \fa$.
The precise expression for $D_{\alpha,\alpha'}$ depends on the rank and is displayed
in \eqref{eq:der} below.    Moreover by \cite{ABG1}, the multiplication in a centerless C$_r$-graded Lie algebra
 $\scl$ is  given by
  \begin{eqnarray}\label{eq:lprod}&&\\
  &&[x\otimes a, y\otimes a']=[x,y]\otimes\frac{1}{2}a\circ a' 
+x\circ y\otimes\frac{1}{2}[a,a']+\tr(xy)D_{a,a'}\nonumber \\ &&[x\otimes a,s\otimes b]
=x\circ s\otimes\frac{1}{2}[a,b]+
[x,s]\otimes\frac{1}{2}a\circ b\nonumber\\ &&[s\otimes b,t\otimes b']=[s,t]\otimes\frac{1}{2}b\circ b'+
s\circ t\otimes\frac{1}{2}[b,b'] +\tr(st)D_{b,b'}\nonumber\\ &&[d,x\otimes a+s\otimes b]
=x\otimes d(a)+s\otimes d(b)\nonumber
\end{eqnarray}
for $x,y\in\g$, $s,t\in\fs$,
$a,a'\in A$, $b,b'\in B$ and $d\in D$,  where

\begin{eqnarray*} w\circ z &=& wz+zw-\frac{1}{r}\tr(wz) \Id \\
{[w,z]}&=&wz-zw \end{eqnarray*}
\noindent for all $w,z \in \mathfrak{gl}_{2r}(\F)$.  Here $\tr$ denotes the usual matrix trace, and
$\Id$ is the identity matrix.   
Note that $x\circ y\in \fs$, $x\circ s\in \g$, $[x,s]\in \fs$,  
$[s,t]\in \g$ and $s\circ t\in \fs$ for $x,y\in\g$, $s,t\in\fs$.
There exists a distinguished element $1\in A$ so that
the grading subalgebra $\g$ of $\scl$ is identified with $\g\otimes 1$,
and
$1\circ \alpha=2\alpha$ and $[1,\alpha]=0$ for all $\alpha\in\fa$.  \medskip

An important remark for the $r=2$ case
is that $s\circ t=0$ for all $s,t\in\fs$,
and so the skew product  on $B$ can be defined arbitrarily in that case. 
This flexibility in defining the skew product is crucial in the determination
of the coordinate algebra in Section 5.  
\medskip

By \cite{ABG1} (see also \cite{Se})  
$\fa$ under the product \eqref{eq:comult}  is an associative algebra for $r \geq 4$
and is an alternative algebra for $r =3$ with $A$  contained in the
nucleus of $\fa$. 
The linear isomorphism $\sigma$,
defined as $a^\sigma=a$ and $b^\sigma=-b$ for $a\in A$ and $b\in B$,   is an involution of $\fa$
since
\begin{eqnarray*} &&A\circ A\subseteq A,\quad [A, A]\subseteq B, \quad A\circ B\subseteq B, \\ &&[A, B]\subseteq
A,\quad  B\circ B\subseteq A,\quad [B,B]\subseteq B.\end{eqnarray*} 
Finally,

\begin{equation}\label{eq:der}
D_{\alpha,\alpha'} = \begin{cases}
\displaystyle{\frac{1}{2r}}\Bigg([L_\alpha,L_{\alpha'}]+ [R_\alpha,R_{\alpha'}]- [L_\alpha,R_{\alpha'}]
+[L_{\alpha^\sigma},L_{{\alpha'}^\sigma}] \\
\hspace{1.4 truein} +[R_{\alpha^\sigma},R_{{\alpha'}^\sigma}]
+[L_{\alpha^\sigma},R_{{\alpha'}^\sigma}]\Bigg)  
&\text{for $r\geq 3$}, \\   \displaystyle{\frac{1}{2}}\Bigg(
[L^+_\alpha,L^+_{\alpha'}]+[L^+_{\alpha^\sigma},L^+_{{\alpha'}^\sigma}]\Bigg)&
\text{for $r= 2$},  \end{cases}
\end{equation}
where $L$ (resp. $R$) is the left (resp. right) multiplication operator on $\fa$,  and $L^+$ is
the left multiplication operator on the plus algebra  $\fa^+$, which is $\fa$  with the multiplication
\begin{equation}\label{eq:cdotmult} \alpha\cdot\alpha':=\hf \alpha\circ \alpha'
\end{equation} for $\alpha,\alpha'\in\fa$.  \medskip

\begin{subsec}{\rm {\textbf{Examples of C$_r$-graded Lie algebras}}} \end{subsec}

Suppose $\fa$ is an algebra with unit element 1 and with product denoted by juxtaposition.
Set  $\alpha \circ \alpha' =  \alpha \alpha' + \alpha' \alpha$
and $[\alpha,\alpha'] = \alpha \alpha' - \alpha' \alpha$  for all $\alpha,\alpha' \in \fa$.
Assume $\fa$ has an involution $\sigma$,  and $A$ (resp. $B$) is the set of 
symmetric (resp. skew-symmetric) elements of $\fa$ relative to $\sigma$.  
Let $\g,\fs$ be as in the previous section, and define $D_{\alpha,\alpha'}$
as in \eqref{eq:der}.     If $\scl = (\g \ot A) \oplus (\fs \ot B) \oplus D_{\fa,\fa}$ under the
multiplication in \eqref{eq:lprod} is a Lie algebra, then we denote it by
$\mathfrak {sp}_{2r}(\fa)$.     In particular, if $\fa$ is any unital associative
algebra with involution having symmetric elements $A$ and skew elements $B$, then $\mathfrak{sp}_{2r}(\fa)$ is a
centerless C$_r$-graded Lie algebra for any $r \geq 2$, and any centerless
C$_r$-graded Lie algebra for $r \geq 4$   is isomorphic to $\mathfrak{sp}_{2r}(\fa)$ for some 
unital associative
algebra $\fa$ with involution.  The centerless C$_3$-graded
Lie algebras are exactly the Lie algebras $\mathfrak{sp}_6(\fa)$, 
where $\fa$ is a  unital alternative algebra with involution whose symmetric
elements $A$ lie in the nucleus of $\fa$.

Now suppose $A$ is a commutative, associative algebra with unit element, and let
$B$ be a left $A$-module.   Assume there is an $A$-bilinear symmetric form
$\zeta: B \times B \rightarrow A$, and define a multiplication on $\mathfrak a = A \oplus B$
by  $(a+b)(a'+b') = aa' + \zeta(b,b') + ab' + a'b$.  Then $\mathfrak a$ with this product
is a Jordan algebra ({\em of Clifford type}).  
For any Jordan algebra $\fa$ of Clifford type, $\mathfrak{sp}_4(\fa)$ 
is a centerless C$_2$-graded Lie algebra.     There is a 
construction described  in \cite{AG} or \cite{Y3}
which results in a B$_r$-graded Lie algebra $ \mathfrak{o}_{2r+1}(\fa)$, and 
when $\fa$ is a Jordan algebra of Clifford type, $\mathfrak{sp}_4(\fa) \cong 
\mathfrak o_5(\fa)$.   \m
 
\begin{subsec} {\rm \textbf{ $\fa^+$ is a Jordan algebra}} \end{subsec}
\m

The plus algebra
$\fa^+=(\fa,\cdot)$   with product
$\alpha\cdot\alpha' =\hf \alpha\circ \alpha'
$ for $\alpha,\alpha'\in\fa$ 
is a Jordan algebra for any alternative algebra 
$\fa$ (see for example, \cite[III, Ex.~3.1]
{M}), hence for the coordinate algebra of any C$_r$-graded Lie algebra, $r \geq 3$ 
By 
\cite[Secs.~2.45 and 2.48]{ABG1} and \cite[Prop.~6.75]{ABG2} (see also \cite[Sec.~4.9]{N3}), the coordinate algebra $\fa$  of any C$_2$-graded Lie algebra
can be identified with the half space $J_{12}$ of
a Jordan algebra $J$ with a {\em triangle}  $(p_1, p_2, q)$. 
Here we show for any Jordan algebra $J$ with a triangle that the half space $J_{12}$ under a suitable
product (see \eqref{eq:2.13}) 
has the structure of a Jordan algebra. 
As a consequence, we obtain that $\fa^+$ is a Jordan algebra for the
coordinate algebra $\fa$ of any  C$_2$-graded Lie algebra.  
 
 \m  Let $J$ be a Jordan algebra with a {\em triangle}  $(p_1, p_2, q)$ 
 and with product denoted by juxtaposition.
 (Facts about triangles quoted here
can be found in \cite[Chap.~III]{J} or \cite[III.6-III.8]{M}.)    
Thus,  the elements  $p_1, p_2, q\in J$ satisfy

\begin{gather} p_1^2=p_1,\quad p_2^2=p_2,\quad  p_1p_2=0,\\
p_1q=\hf q,\quad
p_2q=\hf q,\quad\text{and}\quad q^2=p_1+p_2=1. \end{gather}
We
have the Peirce decomposition  $J=J_{11}\oplus J_{12}\oplus J_{22}$ of $J$,
where $J_{11}, \ J_{22}$,  and $J_{12}$ are the 1, 0,  and $\hf$-eigenspaces, respectively,  of the left multiplication operator  $L_{p_1}$ of the  idempotent $p_1$.  These
spaces  have the following multiplication properties:  

\begin{gather*} J_{11}J_{11}\subseteq
J_{11},\quad  J_{22}J_{22}\subseteq J_{22},\quad  J_{11}J_{22}=0,\\  J_{11}J_{12}+J_{22}J_{12}\subseteq J_{12}\quad
\text{and}\quad J_{12}J_{12}\subseteq J_{11} + J_{22}.\end{gather*}
  Also,
\begin{equation}\label{eq:Peirce} x_{11}(x_{22}x_{12})=x_{22}(x_{11}x_{12})\end{equation}
 for $x_{11}\in J_{11}$, $x_{12}\in J_{12}$,
$x_{22}\in J_{22}$.  The {\em connection involution} determined by the triangle is the transformation
$\sigma:J\rightarrow J$ defined by $\sigma(x)=2(qx)q-x$.  The mapping $\sigma$ is an automorphism of
$J$ of order 2 which stabilizes $J_{12}$, interchanges $J_{11}$ and $J_{22}$, and satisfies $\sigma
L_q=L_q\sigma=L_q$.  We denote its restriction to $J_{12}$ simply by $\sigma$.  Thus,  one can
write $J_{12}=J_{12}^{(+)}\oplus  J_{12}^{(-)}$, where $J_{12}^{(\pm)}$ is the $\pm 1$-eigenspace for
$\sigma$, and $$J_{12}^{(-)}=\{b\in J_{12}\mid qb=0\}.$$  When $J_{12}$ is
the coordinate algebra $\fa = A \oplus B$ of a C$_2$-graded Lie algebra, then
the connection involution  is the involution on
$\fa$ in the previous section, and $A=J_{12}^{(+)}$  and $B=J_{12}^{(-)}$.    Moreover, using the fact that $x_{11}\mapsto x_{11}q$
is a linear isomorphism from $J_{11}$ onto $J_{12}^{(+)}$, we have that the product $\cdot$ on $\fa^+$ is given  by 
\begin{equation}\label{eq:2.13} (a+b)\cdot (a'+b')=\hf \bigg(x_{11}a'+x_{11}'a+x_{11}b'+x_{11}^\sigma
b'+x_{11}'b+(x_{11}')^\sigma b\bigg)-(bb')q\end{equation} for $a, a'\in A$, $b, b'\in B$, $x_{11}, x_{11}'\in
J_{11}$, $a=x_{11}q$ and $a'=x_{11}'q$.

Note that the skew product on $B$ is chosen to be 0 in \cite{ABG1} i.e., $[B,B]=0$, which is essential for
determining the central extensions of a $\hbox{\rm C}_2$-graded Lie algebra, but for any choice of a skew
product on $B$, the plus product is given by the expression in \eqref{eq:2.13}.

\begin{thm}  Let $J=J_{11}\oplus J_{12}\oplus J_{22}$ be a Jordan algebra with a triangle
$(p_1, p_2, q)$. Then the product $\cdot$ on $J_{12}$ defined by \eqref{eq:2.13} coincides with the product of the
$q$-isotope $J^{(q)}$ of $J$ on $J_{12}$,  i.e., $(J_{12}, \cdot)=J_{12}^{(q)}$.  In particular,
$(J_{12}, \cdot)$ is a Jordan algebra. 
\end{thm}

\pf    For $u,v \in J$, the product $\cdot_q$ is defined  by 

$$u \cdot_q v = (uq)v+u(qv)-(uv)q.$$
Thus, if  $a = x_{11}q$, $a' = x_{11}'q  \in J_{12}^{(+)}$ and $b,b' \in J_{12}^{(-)}$, then
\begin{eqnarray*}
&&\hspace{-.3truein}(a+b)\cdot_q (a'+b') \\
&&\hspace{.05truein} =(aq)(a'+b')+(a+b)(qa')-\big((a+b)(a'+b')\big)q\ \  \text{(since $Bq=0$)}\\
&&\hspace{.05truein} =(aq)a'+(aq)b'+a(qa')+b(qa')-(aa')q \\
&& \hspace{.8truein}-(ab')q-(ba')q-(bb')q\\ &&\hspace{.05truein} =(aq)a'+(aq)b'+a(qa')+b(qa')-(aa')q-(bb')q, \end{eqnarray*}
since $(ab')q=(ab')^\sigma q=-(ab')q$ and
$(ba')q=(ba')^\sigma q=-(ba')q$.    Note that
\begin{eqnarray} (aq)a'&=&\big((x_{11}q)q\big)a'=\frac{1}{2}(x_{11}^\sigma+x_{11})a'  \label{eq:calc1} \\
a(qa')&=&a\big(q(x_{11}'q)\big)=\frac{1}{2}a\big((x'_{11})^\sigma+x_{11}'\big),  
\end{eqnarray}  but $x_{11}^\sigma a'=x_{11}'a$
and $(x'_{11})^\sigma a=x_{11}a'$ by \eqref{eq:Peirce}, and hence $(aq)a'+a(qa')=x_{11}a'+x_{11}'a$.  Also, we have
$(aq)b'=\big((x_{11}q)q\big)b'=\frac{1}{2}(x_{11}^\sigma+x_{11})b'$ and
$b(qa')=b\big(q(x_{11}'q)\big)=\frac{1}{2}b\big((x'_{11})^\sigma+x_{11}'\big)$.  
Thus, it is enough to show that
\begin{equation}\label{eq:calc2} (aa')q=\frac{1}{2}\big((aq)a'+a(qa')\big). \end{equation}

We use the following two identities to establish
\eqref{eq:calc2}:  
\begin{eqnarray} (qx_{ii})x_{12}&=&(qx_{12})x_{ii}+(q(x_{12}x_{ii}))p_j  \label{eq:calc3} \\
(x_{ii}y_{ii})x_{12}&=&(x_{ii}x_{12})y_{ii} \label{eq:calc4}  \end{eqnarray}  for $x_{ii}, y_{ii}\in J_{ii}$,
$i,j\in\{1,2\}$, $i\neq j$, and $x_{12}\in J_{12}$.

Now for \eqref{eq:calc3}, the formula \cite [(1.3.3)]{MN}, which was stated for
Jordan triple systems, can be adapted for Jordan algebras to say
\begin{eqnarray*} && \hspace{-.45truein} (mx_{12})x_{ii}+m(x_{12}x_{ii})-(mx_{ii})x_{12}\\
&&\hspace{.8truein}= (m(x_{ii} x_{12}))p_i+m((x_{ii}
x_{12})p_i)-(mp_i)(x_{ii} x_{12})\end{eqnarray*}  for $m\in J_{12}$.  Let $m=q$.  
Then, since
$p_i(x_{ii}x_{12})=\hf x_{ii}x_{12}$ and $qp_i=\hf q$, we have 
\begin{eqnarray*}
&& \hspace{-.45truein}(qx_{12})x_{ii}+q(x_{12}x_{ii})-(qx_i)x_{12}\\
&& \hspace{.8truein} =(q(x_{ii}x_{12}))p_i+\frac{1}{2}q(x_{ii}
x_{12})-\frac{1}{2}q(x_{ii} x_{12}) \\ && \hspace{.8truein} = (q(x_{ii} x_{12}))p_i.\end{eqnarray*}  
Note that $q(x_{ii} x_{12})\in
J_{12} J_{12}\subset J_{11}\oplus J_{22}$, and so $(q(x_{ii} x_{12}))=(q(x_{ii} x_{12}))p_1+(q(x_{ii}
x_{12}))p_2$. Hence, $(qx_{12})x_{ii}+(q(x_{12}x_{ii}))p_j-(qx_{ii})x_{12}=0$, which is \eqref{eq:calc3}.  Also observe that 
$(qx_{12})x_{ii}=((qx_{12})x_{ii})p_i$ since $qx_{12}\in J_1\oplus J_2$ and
$J_1J_2=0$.  Thus, \eqref{eq:calc4} follows from applying
\cite [(1.3.5)]{MN} to Jordan algebras.
 
In demonstrating \eqref{eq:calc2}, we write $y=y_1+y_2$  for $y\in J_{11}\oplus J_{22}$     
$(y_i \in J_{ii})$  to simplify
the notation.     Then we have 
\begin{eqnarray*}  (aa')q &=&((qx_{11})a')q\\ &=&((qa')x_{11})q+(q(a'x_{11}))_2q\quad \ \text{by \eqref{eq:calc3}}\\
&=&((qa')_1x_{11})q+(q(a'x_{11}))_2q\quad \text{since $J_{22}J_{11}=0$}\\
&=&(qa')_1(x_{11}q)+(q(a'x_{11}))_2q\quad \text{by \eqref{eq:calc4}}\\ &=&(qa')_1a+(q(a'x_{11}))_2q.  
\end{eqnarray*}
Note that $a=qx_{11}=qx_{11}^\sigma$ and $x_{11}^\sigma\in J_{22}$, and so, by a similar argument, we
also get
$(aa')q=(qa')_2a+(q(a'x_{11}^\sigma))_1q$. Hence,
$$2(aa')q=(qa')a+(q(a'x_{11}))_2q+(q(a'x_{11}^\sigma))_1q.$$  Since $L_q=L_q\sigma$, we have
$q(a'x_{11}^\sigma)=q(a'x_{11})$.  Hence,  by  \eqref{eq:calc1}, 
$$(q(a'x_{11}))_2q+(q(a'x_{11}^\sigma))_1q=(q(a'x_{11}))q
=\frac{1}{2}(a'x_{11}^\sigma+a'x_{11})=a'(qa).$$  Thus,  \eqref{eq:calc2}  holds, and the proof is 
finished. 
\qed   \m

 Suppose that $J$ is a Jordan algebra with a triangle $(p_1,p_2,q)$
and connection involution $\sigma$.
Let  $\fa=(J_{12}, \cdot)$, $A=J_{12}^{(+)}$,   and $B=J_{12}^{(-)}$, where the product
$\cdot$ is as in \eqref{eq:2.13}.  Define a new multiplication  
on $\fa$ by
\begin{eqnarray*} a a' &=& x_{11} \cdot a', \\
a b &=& x_{11}\cdot b, \\
b a &=& x_{11}^\sigma\cdot b, \\
bb' &=& -(b\cdot b')\cdot q \end{eqnarray*}
for $a = x_{11}\cdot q$, $a' \in A$, and $b,b' \in B$ as in \cite[Sec.~2.48]{ABG1} so that
$\alpha \cdot \alpha' = \hf(\alpha\alpha' + \alpha' \alpha)$ for all $\alpha,\alpha' \in \fa$.
Thus, $\fa$ is the coordinate algebra of a Lie algebra graded by  $\hbox{\rm C}_2$,
and every coordinate algebra $\fa$  of a Lie algebra graded by $\hbox{\rm C}_2$ has this form,
(see the discussion in \cite[Sec.~2.51]{ABG1}).   So in summary, we have

\begin{cor}\label{cor:2.20} Let $\fa = A \oplus B$ be a coordinate algebra of a Lie algebra graded by  $\hbox{\rm C}_2$.
Then  $\fa^+ = (\fa,\cdot)$ is a Jordan algebra with involution.  In particular, $\fa$
is a Jordan admissible algebra with involution (for any choice of skew product on $B$).
\end{cor}

\begin{rem}\label{rem:2.21}{\rm In \cite{Se},  Seligman  proved that $(A,\cdot)$ is a Jordan algebra by the following argument. 
(Actually Seligman was working with finite-dimensional simple Lie algebras
graded by C$_2$, but the same proof applies in the general setting.)   The
universal relation 
$\tr(x[y,z])=\tr(y[z,x])=\tr(z[x,y])$ implies that 
$$ D_{a,a'\circ a''}+D_{a',a''\circ a'}+D_{a'',a\circ a'}=0
$$ for $a,a',a''\in A$. In particular,
$D_{a,a^2}=0$.    Moreover,  $$D_{a,a'}a''=\frac{2}{r}\Big(a\circ(a'\circ a'')-a'\circ(a\circ
a'')\Big).$$   Combining those results gives
$a\cdot(a^2\cdot a')=a^2 \cdot(a \cdot a')$ for all $a,a'\in A$ so that $(A,\cdot)$ is a Jordan algebra.
In our classification of Lie $G$-tori of type $\hbox{\rm C}_2$ in Section 5, we only require the fact that $(A,\cdot)$ is a Jordan algebra.   However,
we have included the proof that $(J_{12}, \cdot)$ and $\fa^+ =(\fa,\cdot)$ are Jordan algebras, 
as those results may be of independent interest.}
\end{rem}  \bigskip

\section {The coordinate algebra of a $(\hbox{\rm C}_r, G)$-graded Lie algebra}
\m

In this section,  we show that the coordinate algebra of a $(\hbox{\rm C}_r, G)$-graded Lie algebra $\scl$ is $G$-graded.   More specifically, we prove the following theorem.
\m

\begin{thm}\label{Thm:3.5} \begin{itemize}
\item[{\rm (i)}]  Let $\scl$  be a  $(\hbox{\rm C}_r, G)$-graded Lie algebra.
Then for $r \geq 3$, the  coordinate algebra $\fa=A\oplus B$ of $\scl$  is $G$-graded 
and has a graded involution.  
Also, for $r \geq 2$,  \ $(A,\cdot)$ is an  $\langle L \rangle$-graded Jordan algebra 
with graded involution, where $\langle L \rangle$ is the subgroup
of $G$ generated by  $L = \{g \in G \mid \scl_\mu^g \neq 0, \mu \in \Delta_{lg}\}$.
\item[{\rm (ii)}]  $\mathfrak{sp}_{2r}(\fa)$ is  a centerless $(\hbox{\rm C}_r, G)$-graded 
Lie algebra  for any $G$-graded associative algebra $\fa$  with graded involution 
if $r\geq 2$, or for any $G$-graded alternative algebra $\fa$  with graded involution whose  
symmetric elements $A$ are in the nucleus of $\fa$ if $r=3$.
In addition, $\mathfrak{sp}_{4}(\fa)$ is  a centerless $(\hbox{\rm C}_2, G)$-graded 
Lie algebra  for any $G$-graded Jordan algebra $\fa$  of Clifford type.  \end{itemize}
\end{thm}

\pf   (i) We suppose that  $\scl$ is a  $(\hbox{\rm C}_r, G)$-graded Lie algebra. 
Thus, we are assuming that $\scl$ is $\Delta$-graded, 
$\scl= \bigoplus_{\mu \in \Delta \cup \{0\}} \scl_\mu$  
for $\Delta =$ C$_r$  $(r \geq 2)$ with
grading subalgebra  $\g$;   \  $\scl$ is $G$-graded  $\scl = \bigoplus_{g \in G} \scl^g$ and has a decomposition 
$$\scl=\bigoplus_{\mu\in\Delta\cup\{0\}}\ \bigoplus_{g\in G}\ \scl_\mu^g,$$ 
where $\scl_\mu^g = \scl_\mu \cap \scl^g$; \  $\g \subseteq \scl^0$; \ and
$\supp \scl$ generates $G$.    Then  we have
\begin{equation}\label{eq:3.1}
\scl_\mu=
\begin{cases}
\g_\mu\otimes A&\text{if $\mu\in\Delta_{lg}$}\\ (\g_\mu\otimes A)\oplus (\fs_\mu\otimes B)
&\text{if $\mu\in\Delta_{sh}$}.
\end{cases} \end{equation}  \smallskip

Set $$L=\{g\in G\mid \scl_\mu^g\neq 0,\mu\in\Delta_{lg}\}$$
and let $\langle L \rangle$ be the subgroup of $G$ generated
by $L$.    For all $\mu\in\Delta_{lg}$ and $g\in G$,
we define $A_\mu^g$ using
$$
\scl_\mu^g=\g_\mu\otimes A_\mu^g.
$$ Then 
$A=\bigoplus_{g\in G}\ A_\mu^g$, and in particular, $A_\mu^g=0$ if $g\notin L$. For any
$\mu,\nu\in\Delta_{lg}$, there exist $\gamma_1,\gamma_2\in\Delta_{sh}$  such that
$\nu+\gamma_1\in\Delta_{sh}$ and
$\mu=\nu+\gamma_1+\gamma_2$. Then,
$$
\g_\mu\otimes A_\mu^g =[[\g_\nu\otimes A_\nu^g,\g_{\gamma_1}\otimes 1],
\g_{\gamma_2}\otimes 1] =
\g_\mu\otimes\ A_\nu^g
\quad\text{for all $g\in G$}.
$$ Therefore,
$A_\mu^g=A_\nu^g$ for all $\mu,\nu \in \Delta_{lg}$, and for $g \in G$,  we specify that 
$$ A^g:=A_\mu^g
\quad\text{for any choice of $\mu\in\Delta_{lg}$.}
$$
Then 
\begin{equation}\label{eq:3.2} A=\bigoplus_{g\in G}\ A^g=\bigoplus_{l\in L}\ A^l, \end{equation}
where 
$\scl_\mu^l=\g_\mu\otimes A^l$
for all $\mu\in\Delta_{lg}$ and $l\in L$.  
The algebra $A$ is graded by the group $\langle L\rangle$,  and $1\in A^0$.  \m 
 
Let $\mu\in\Delta_{sh}$ and $g\in G$. We define $A_\mu^g$ and $B_\mu^g$ via the relations
$$\g_\mu\otimes A_\mu^g=(\g_\mu\otimes A)\cap\scl_\mu^g
\quad\text{and}\quad
\fs_{\mu}\otimes B_\mu^g=(\fs_{\mu}\otimes B)\cap\scl_\mu^g.
$$ We claim that
\begin{equation}\label{eq:claim}
\scl_\mu^g= (\g_\mu\otimes A_\mu^g)\oplus (\fs_\mu\otimes B_\mu^g)
\quad\text{and}\quad A_\mu^g=A^g.
\end{equation} 
To see this, let $w\in \scl_\mu^g$.    Then by \eqref{eq:3.1},
$w=u + v$ for some $u  \in\g_\mu\otimes A$ and
$v\in \fs_\mu\otimes B$.  We need to show that $u,v \in \scl_\mu^g$. Now by \eqref{eq:3.2}, we have
$u=\sum_{l\in L}\ e_\mu\otimes a_l$ for  some $0\neq e_\mu\in \g_\mu$ and
$a_l\in A^l$.     We can find some 
$\nu\in\Delta$ such that $\mu+\nu\in\Delta_{lg}$.    Let 
$0\neq e_\nu\in\g_\nu=\g_\nu\otimes 1\subset \scl_\nu^0$. Then
\begin{equation}\label{eq:3.3} [e_\nu,w]\in\scl_{\mu+\nu}^g=
\begin{cases} 0&\text{if $g\notin L$}\\
\g_{\mu+\nu}\otimes A^g&\text{if $g\in L$}.
\end{cases} \end{equation}
 Let $0\neq s_\mu\in\fs_\mu$.  Since weight spaces of $\fs$ relative to 
 $\mathfrak h$  are one-dimensional,  $v=s_\mu\otimes b$ for some $b\in B$, and
$$[e_\nu,v]=[e_\nu,s_\mu\otimes b] =[e_\nu, s_\mu]\otimes b\ \Big(+e_\nu\circ s_\mu\otimes
\hf[1,b]\Big).$$  But
$[e_\nu, s_\mu]\otimes b=0$ since $\mu+\nu\in\Delta_{lg}$, and hence $[e_\nu,v]=0$. Thus we obtain
$[e_\nu,u]=[e_\nu,w]$.

If $g\notin L$, then, by \eqref{eq:3.3},
$0=[e_\nu,w]=[e_\nu,u]=\sum_{l\in L}\ [e_\nu,e_\mu\otimes a_l]=0$, and so $[e_\nu,e_\mu\otimes a_l]=0$
for all $l\in L$. Since $[\g_\nu,\g_\mu]\neq 0$, we get $a_l=0$ for all $l\in L$, i.e., $u=0$.
Therefore, $w=v\in \scl_\mu^g$.

If $g\in L$, then
$[e_\nu,u]=[e_\nu,w]\in\scl_{\mu+\nu}^g$, and so $u=e_\mu\otimes a_g\in\g_\mu\otimes A^g$. Note
that there exists $\gamma \in\Delta_{lg}$ such that
$\mu-\gamma\in\Delta$. So by \eqref{eq:3.2}, we have
$$
\g_\mu\otimes A^g =[\g_\gamma\otimes A^g,\g_{\mu-\gamma}
\otimes 1]
\subseteq [\scl_\gamma^g,\scl_{\mu-\gamma}^0]
\subseteq  \scl_\mu^g.
$$ Therefore, $u \in \scl_\mu^g$ and $v=w-u \in \scl_\mu^g$. Finally, since $\mu+\nu\in\Delta_{lg}$, it
follows that
$$ [\g_\mu\otimes A_\mu^g,\g_{\nu}
\otimes 1] =
\g_{\mu+\nu}\otimes A_\mu^g
\subseteq
\g_{\mu+\nu}\otimes A^g.
$$ Hence $A_\mu^g\subseteq A^g$. Also,
$$ [\g_{\mu+\nu}\otimes A^g,\g_{-\nu}
\otimes 1] =
\g_{\mu}\otimes A^g
\subseteq
\g_{\mu}\otimes A_\mu^g,
$$ and so $A^g\subseteq A_\mu^g$. Thus our claim \eqref{eq:claim}  is settled.

Now,
$B=\bigoplus_{g\in G}\ B_\mu^g$, and
$\fs_\mu\otimes B_\mu^g=
\scl_\mu^g$ if $g\notin L$ since $A^g=0$. If $\mu,\nu\in\Delta_{sh}$ and
$\mu-\nu\in\Delta$, then
$$
\fs_\mu\otimes B_\mu^g =[\fs_\nu\otimes B_\nu^g,\g_{\mu-\nu}
\otimes 1] =
\fs_\mu\otimes B_\nu^g
\quad\text{for all $g\in G$}.
$$ Therefore,
$B_\mu^g=B_\nu^g$. Thus by the same argument as in \cite[(5.11)]{AG},  we get $B_\mu^g= B_\nu^g$ for any
$\mu,\nu\in\Delta_{sh}$ and all $g\in G$. So for $g\in G$ we put
$$ B^g:=B_\mu^g
\quad\text{for any choice of $\mu\in\Delta_{sh}$.}
$$ Consequently, 
$$ B=\bigoplus_{g\in G}\ B^g,
\quad\text{with}
$$
$$
\scl_\mu^g=
\fs_\mu\otimes B^g
\quad\text{for all $\mu\in\Delta_{sh}$ and $g\notin L$},
$$ and
\begin{equation}\label{eq:3.4}
\fs_\mu^g= (\g_\mu\otimes A^g)
\oplus (\fs_\mu\otimes B^g)
\quad\text{for all $\mu\in\Delta_{sh}$ and $g\in L$}.\end{equation}

Let
$$
\fa=\bigoplus_{g\in G}\ \fa^g,
\quad\text{where}\quad
\fa^g:= A^g\oplus B^g
$$ ($A^g=0$ if $g\notin L$). Let $S:=\supp\fa=\supp\scl_\mu$ for $\mu\in\Delta_{sh}$.   By \eqref{eq:3.4}, we
have $L\subseteq S$, and so
$S+S\supset\supp\scl$, which generates $G$. Hence $S$ generates $G$.

At this stage we know that $\fa$ is a vector space graded by the group $G$,
and the support of $\fa$ generates $G$.  
We need to verify that $\fa$ is a graded algebra.  
Now when $w = E_{1,2} -E_{2r-1,2r}\in \g_{\varepsilon_1-\varepsilon_2}$
and $ z = E_{2,2r-1} \in \g_{2 \varepsilon_2}$, 
we have $[w,z] =E_{1,2r-1}+ E_{2,2r} \in \g_{\varepsilon_1+\varepsilon_2}$ and
 $w \circ z =   E_{1,2r-1}- E_{2,2r}
 \in \fs_{\varepsilon_1+\varepsilon_2}$, both of which are nonzero. 
  Then the product
 $[w \otimes a, z \otimes a']$ with $a\in A^g$, $a' \in A^h$ shows that $A^g \circ  A^h \subseteq A^{g+h}$ 
and $[A^g, A^h]  \subseteq B^{g+h}$, which combine to say $A^g A^h \subseteq \fa^{g+h}$.  
 
    The elements $s = E_{2,1}+E_{2r,2r-1}$ and $s' = 
E_{1,2r-1}-E_{2,2r}$
 belong to $\fs$  as does $t = 
 E_{1,3} + E_{2r-2,2r}$ when $r\geq 3$.   
 Setting $x = E_{1,2} -E_{2r-1,2r} \in \g$, we have
 $[x,s] = E_{1,1}-E_{2,2} -E_{2r-1,2r-1} + E_{2r, 2r}$ and
 $x \circ s =  E_{1,1}+E_{2,2} -E_{2r-1,2r-1} - E_{2r,2r}$, from which
 we can deduce that $A^g \circ B^h \subseteq B^{g+h}$ and $[A^g,B^h] \subseteq A^{g+h}$.
 Thus, $A^g B^h \subseteq \fa^{g+h}$.   We can use the fact that
 $[s,s'] = 2E_{2,2r} \in \g_{2\varepsilon_2}$ to determine
 that $B^g \circ B^h \subseteq A^{g+h}$.  
 Now when $r \geq 3$,   $s \circ t = E_{2,3}+ E_{2r-2,2r-1} \neq 0$, 
 from which we obtain  $[B^g,B^h] \subseteq B^{g+h}$.
 Thus, for $r \geq 3$, we have   $B^g B^h \subseteq \fa^{g+h}$.
 The product $s \circ t$ is identically 0 on $\fs$ when $r =2$,  and all
 we can deduce in this case is that $B^g \circ B^h \subseteq A^{g+h}$.  
 
 These arguments have shown that $\fa$ is graded for $r\geq 3$.   By Corollary \ref{cor:2.20} or Remark
\ref{rem:2.21},  it follows that $(A,\cdot)$ is a Jordan $\langle L\rangle$-graded algebra for $r \geq 2$.
The involution $\sigma$ is clearly graded,  so we have (i).  

(ii) For this second part, 
let $\scl=\mathfrak{sp}_{2r}(\fa)$, where
$\fa=A\oplus B=\bigoplus_{g\in G}\ (A^g\oplus B^g)$  is a $G$-graded algebra with
symmetric elements $A$  and   skew elements $B$ relative to a graded involution.  
Assume $\fa$ is associative;   or in the $r = 3$ case,   an
alternative algebra such that  $A$ lies in the nucleus of
$\fa$;  or in the case $r = 2$, a Jordan algebra.    For $g\in G$, set
$\scl_\mu^g:= (\g_\mu\otimes A^g)\oplus (\fs_\mu\otimes B^g)$ if $\mu\in\Delta_{sh}$, and
$\scl_\mu^g:=
\g_\mu\otimes A^g$ if $\mu\in\Delta_{lg}$. Then
$\scl$ admits a compatible $G$-grading,  $\scl = \bigoplus_{g \in G} \scl^g$,  with 
$\scl^g =\bigoplus_{\mu\in\Delta\cup\{0\}} \scl_\mu^g$ and  $\scl_0^g=\sum_{\mu\in\Delta}\ \sum_{g=g'+g''}\ [\scl_\mu^{g'},\scl_{-\mu}^{g''}]$, so that 
$\scl$ is a 
$(\hbox{\rm C}_r,G)$-graded Lie algebra.\qed 
\m

\section{The coordinate algebra of \\
 a division (C$_r, G$)-graded Lie algebra} 

In this section we investigate Lie algebras that are $(\hbox{\rm C}_r,G)$-graded
and satisfy the division property (see Remarks \ref{rems:2.2});   that is, 
the so-called {\it division $(\hbox{\rm C}_r,G)$-graded Lie algebras}.    Our main
result will be that the coordinate algebra of such a Lie algebra  is a  division $G$-graded algebra where by that we mean the following.

\begin{defn}\label{def:Gtorus} A $G$-graded unital (associative, alternative, or Jordan) algebra $\mathcal A$  is
said to be {\em division $G$-graded} (or have the division property) if all nonzero homogeneous
elements are invertible.  If $\mathcal A$ is  a division $G$-graded algebra such 
that $\dim \mathcal A^g \leq 1$ for all homogeneous
spaces $\mathcal A^g$, then $\mathcal A$ is a {\em $G$-torus}.
A $\mathbb Z^n$-torus is referred to as an {\em $n$-torus} or simply a {\em torus}.
\end{defn} 

\begin{exams} {\rm (1) An associative $G$-torus is nothing but a twisted group algebra
$\F^t[G]$.    Thus, an associative $n$-torus (also known in the literature as a {\em quantum torus})  is a  Laurent polynomial ring $\F_{\underline q}[t_1^{\pm 1}, \dots, t_n^{\pm 1}]$
in $n$ variables with multiplication given  by $t_i t_j = q_{i,j} t_j t_i$
where $\underline q = (q_{i,j})$ is an $n \times n$ matrix with entries in $\F^\times$ such that
$q_{i,i} = 1$ for all $i$ and $q_{j,i}= q_{i,j}^{-1}$. \smallskip

(2)  Alternative tori were classified in \cite{BGKN} for fields $\F$
such that every element of $\F$ has a square root in $\F$,  and in  \cite{Y2} for arbitrary  $\F$.
An alternative torus is either a quantum torus or an octonion torus (sometimes
called a Cayley torus).   In the second case, 
it is the Cayley-Dickson algebra
$\mathbb O_n$
over the ring of Laurent polynomials 
$\F[t_1^{\pm 1}, \ldots,t_n^{\pm 1}]$ in $n$ variables 
for some $n \geq 3$ obtained by successively
applying the Cayley-Dickson process with elements
$x_1,x_2,x_3$, such that $x_i^2 = t_i$, where the $t_i$  
are the structure constants of the process.
A  graded involution $\sigma$ of $\mathbb O_n$
whose symmetric elements lie in the nucleus,
must be the standard involution, i.e.,
the involution determined by $x_i\mapsto -x_i$ for $1 \leq i \leq 3$, 
and $t_i\mapsto t_i$ for $4\leq i\leq n$. \m

(3) If a division $G$-graded Jordan algebra 
$\fa=\bigoplus_{g\in G}\ \fa^g$
has a decomposition $\fa^g=A^g\oplus B^g$ for each $g\in G$
so that 
$A=\bigoplus_{g\in G}\ A^g$ is a commutative  associative subalgebra, and
$\fa=A\oplus B$
is a Jordan algebra over $A$ of a symmetric bilinear form
on the graded $A$-module $B:=\bigoplus_{g\in G}\ B^g$, then
we say that $\fa$ is a {\em division $G$-graded Jordan algebra of Clifford type}.
The algebra $\fa$  has a natural involution $\sigma$,  with
$\sigma(a)=a$ for all $a\in A$ and $\sigma(b)=-b$ for all $b\in B$.
We call this $\sigma$ the standard involution.
If $\fa$ is a Jordan $G$-torus, then it is said to be 
a {\it Clifford $G$-torus}.}
\end{exams}

Suppose that $\scl$ is a division $(\hbox{\rm C}_r,G)$-graded Lie algebra with
grading subalgebra $\g$.  
We may assume that $\scl$ is centerless,
as the coordinate algebras of $\scl$ and
$\scl/\mathcal Z(\scl)$ are the same  because
the center $\mathcal Z(\scl)$ is contained in the sum of
the trivial $\g$-submodules of $\scl$.    By the previous section,
the coordinate algebra $\fa = A \oplus B$ of $\scl$ is $G$-graded.  

 Let $0\neq a+b\in A^g\oplus B^g=\fa^g$ for 
$a\in A^g$,
$b\in B^g$ and
$g\in S = \hbox{\rm supp}\,\fa$ (Note $a=0$ if $g\notin L$).  Let $\mu:=\varepsilon_1-\varepsilon_2\in\Delta_{sh}$, 
$e:=E_{1,2}-E_{3,4}$, $e':=\hf(E_{2,1}-E_{4,3})$
$s:=E_{1,2}+E_{3,4}$ and $s':=\hf(E_{2,1}+E_{4,3})$.
Then for $r \geq 2$, we have

\begin{eqnarray*} 
{[e,e']} &=&[s,s'] ={\hf }\mu^\vee  \quad (\text{recall we are assuming $\mathcal Z(\scl) = 0$}), \\ 
e\circ s' &=& s\circ e', \qquad \quad 
\ \text{(which is linearly independent of $\mu^\vee$),}\\ 
{[e, s']}&=&[s, e']\neq 0,  \quad 
\ \text{and}\\
\tr(ee')&=&\tr(ss')\neq 0.
\end{eqnarray*}

Also, one can check that if $r \geq 3$, then
$$ e \circ e'  =s \circ s' ,
\quad (\text{which is linearly independent of $[e , s' ]$}), 
$$ and 
$e \circ e'  =s \circ s' =0$ if $r = 2$.

Now, $e_\mu\otimes a+s_\mu\otimes b\in\scl_{\mu}^g$, and by the division property of $\scl$, there
exists 
$y\in\scl_{-\mu}^{-g}$ such that $[e_\mu\otimes a+s_\mu\otimes b,y]=\mu^\vee$. Since 
$\scl_{-\mu}^{-g}= (\g_{-\mu}\otimes A^{-g})\oplus (s' \otimes B^{-g})$, we have that
$y=e' \otimes a'+s' \otimes b'$ for suitable elements  $a'\in A^{-g}$ and $b'\in B^{-g}$.  Consequently,

\begin{eqnarray*} 
\mu^\vee\otimes 1=\mu^\vee &=&[e \otimes a+s \otimes b,e' \otimes a'+s' \otimes b']\\
&=&[e ,e' ]\otimes \frac{1}{2}a\circ a'+ e \circ e' \otimes \frac{1}{2}[a, a']+
\tr(e e' )D_{a,a'}\\ &&+e \circ s' \otimes \frac{1}{2}[a,b']+ [e ,
s' ]\otimes \frac{1}{2}a\circ b'\\ &&+s \circ e' \otimes \frac{1}{2}[b,a']+ [s ,
e' ]\otimes \frac{1}{2}b\circ a'\\ &&+[s ,s' ]\otimes \frac{1}{2}b\circ b'+ s \circ
s' \otimes \frac{1}{2}[b, b']+
\tr(s s' )D_{b,b'}, \end{eqnarray*}
from which we deduce that 
\begin{eqnarray*} &&a\circ a'+ b\circ b'=2,\\ &&[a,b']+
 [b,a']=0=a\circ b' +b\circ a',\\ &&D_{a,a'}+D_{b,b'}=0, \quad \text{and} \\ &&  
[a, a']+[b, b']=0   \quad (\text{if $r \geq 3$})
\end{eqnarray*}
 So we get $(a+b)\circ (a'+b')=2$ and
$D_{a+b,a'+b'}=D_{a,a'}+D_{b,b'}=0$. Moreover,  if $r\geq 3$, then
$[a+b, a'+b']=0$. Therefore, $a+b$ is invertible in the Jordan algebra $\fa^+$
 if $r\geq 2$, and $a+b$ is also invertible in the alternative algebra $\fa$ if $r\geq 3$. Thus,

 \begin{equation} \label{eq:items4}   \end{equation}
 \begin{itemize} 
 \item[{\rm (i)}]
$\fa^+$ is a division $G$-graded Jordan algebra with graded involution for $r\geq 2$.
 \item[{\rm (ii)}] $(A,\cdot)$ is a division $\langle L\rangle$-graded Jordan algebra for $r\geq 2$,
 \item[{\rm (iii)}] $\fa$ is a division $G$-graded alternative algebra with graded involution if $r\geq 3$, and 
  \item[{\rm (iv)}]$\fa$ is a division $G$-graded associative algebra with graded involution if $r\geq 4$.
\end{itemize}  
  
   In particular,
$S=\supp\scl_\mu$ for $\mu\in\Delta_{sh}$ and
$L=\supp\scl_\nu$ for $\nu\in\Delta_{lg}$ are reflection spaces of $G$ (in the sense of \cite{Y3}), and
$S=G$ if $r\geq 3$. Thus we have established the following result.  \bigskip

\begin{thm}\label{Thm:4.2} Let $\scl$ be a centerless division $(\hbox{\rm C}_r, G)$-graded Lie algebra. Then  $\scl\cong\mathfrak{sp}_{2r}(\fa)$ for some division $G$-graded algebra $\fa$ with
graded involution  such that (i)-(iv) of \eqref{eq:items4} hold.  Also, $\mathfrak {sp}_{2r}(\fa)$ is  a centerless division $( 
\hbox{\rm C}_r, G)$-graded Lie algebra  for any division $G$-graded associative algebra $\fa$ with graded involution if
$r\geq 2$, for any division $G$-graded alternative algebra $\fa$ with graded involution so that the
symmetric elements are in the nucleus of $\fa$ if $r=3$,
or for any division $G$-graded Jordan algebra $\fa$ of Clifford type if $r=2$.
\end{thm}

\pf    All this is apparent from our discussions above, except perhaps for the division
property of $\scl = \mathfrak {sp}_{2r}(\fa)$ in the second statement.
For $\mu\in\Delta_{lg}$ and $g\in L$, let $e\in\g_\mu$ and
$e'\in \g_{-\mu}$ be such that
$[e,e']=\mu^\vee.$  Then for $0\neq v\in\scl_\mu^g$, there exists $0\neq a\in A^g$ such that
$v=e\otimes a$. Taking $w=e'\otimes a^{-1}\in\scl_{-\mu}^{-g}$, we get 
$[v,w]=\mu^\vee$. (Note that $[a,a^{-1}]=aa^{-1}-a^{-1}a=0$.)

For $\mu\in\Delta_{sh}$, it is easy to see the existence of
the elements $e \in\g_{\mu}$,
$e' \in \g_{-\mu}$,
$s \in\fs_{\mu}$ and $s' \in \fs_{-\mu} $ satisfying (4.1). 
Then for $g\in S$ and $0\neq v\in\scl_\mu^g$,
there exist $a\in A^g$ and $b\in B^g$ such that
$v=e \otimes a+s \otimes b$. Taking $w=e' \otimes a'+s' \otimes b'\in\scl_{-\mu}^{-g}$,
where 
$(a+b)^{-1}=a'+b'$, we get 
$[v,w]=\mu^\vee$. Hence $\scl$ is division graded.
\qed  \bigskip 

\begin{rem}\label{rem:4.3} {\rm 
The argument above affords a more direct and easier proof  that the division property holds for the coordinate algebra of
a centerless division $(\hbox{\rm C}_r, G)$-graded Lie algebra  
than the one given in \cite[pp.~99-101]{Se}, which treats 
only a particular case of this result; namely, 
 that the coordinate algebra of
a finite-dimensional simple Lie algebra  of relative type C$_r$
is a division algebra.}
\end{rem} \medskip

\section{The coordinate algebra of a Lie $G$-torus of type $\hbox{\rm C}_2$}  \medskip

We apply the following identities 
to determine the coordinate algebra $\fa = A \oplus B$  of a Lie $G$-torus of type $\hbox{\rm C}_2$. 
Seligman \cite[pp.~88-95]{Se} used these same identities in his
classification of the finite-dimensional simple Lie algebras of characteristic zero
graded by  the root system $\hbox{\rm C}_2$, but they are valid in
any  $\hbox{\rm C}_2$-graded Lie algebra.   In expressing them, 
we have translated them into our notation using
 $\circ$ and $[\cdot,\cdot]$ and have written the inner derivations as left
 operators rather than right operators as in \cite{Se}.  
Each identity carries two numbers - the left being the reference in \cite{Se}
and the right being our own equation label.
 
\begin{eqnarray}
&&\hspace{-.6truein}(37')\quad \ a\circ (a''\circ a')-a''\circ(a\circ a') =[a,[a'',a']]-[a'',[a,a']],  \label{eq:5.1}\\
&&\hspace{-.6truein}(38'')\quad  [a,a'\circ a'']=[a,a']\circ a''-[a'',a]\circ a',  \label{eq:5.2}\\
&&\hspace{-.6truein}(39)\quad \ \  {[[a',a''],a]}=a'\circ (a''\circ a)-a''\circ (a\circ a')\label{eq:5.3} \\
&&\hspace{-.6truein}(39')\hspace{.9truein} =4D_{a',a''},  \nonumber  \\ 
&&\hspace{-.6truein}(40)\ \ \ D_{[a,a'],b}=D_{[b,a'],a}-D_{[b,a],a'}, \label{eq:5.4}\\
&&\hspace{-.6truein}(41'')\  [b,a\circ a']=[b,a]\circ a'- [a',b]\circ a, \label{eq:5.5}\\
&&\hspace{-.6truein}(42')\ \  [b\circ a, a']
=[b,a\circ a']+b\circ[a,a']-[b,a]\circ a',\label{eq:5.6} \\
&&\hspace{-.6truein}(42'')\ \  a'\circ[b,a]=b\circ[a,a']+[b\circ a',a], \label{eq:5.7} \\ 
&&\hspace{-.6truein}(43')\ \ 4D_{a,a'}b=[a,[a',b]]-[a',[a,b]], \label{eq:5.8} \\
&&\hspace{-.6truein}(44')\ \  [a,[b,a']]=(b\circ a)\circ a'-b\circ (a\circ a'), \label{eq:5.9} \\
&&\hspace{-.6truein}(46)\ \ \ D_{b\circ b',a}+D_{b'\circ a,b}+D_{b\circ a,b'}=0, \label{eq:5.10} \\
&&\hspace{-.6truein}(51)\ \ \ b\circ (b'\circ a) -b'\circ(b\circ a)=4D_{b,b'}a=[b,[b',a]]-[b',[b,a]],\label{eq:5.11}\\
&&\hspace{-.6truein}(52)\ \ \ 2a\circ(b\circ b')= b\circ(b'\circ a)+b'\circ(b\circ a) +[b,[b',a]]+[b',[b,a]],
\label{eq:5.12}\\
&&\hspace{-.6truein}(53)\ \ \ [a,b\circ b']=[a,b]\circ b'-[b',a]\circ b,\label{eq:5.13} \\
&&\hspace{-.6truein}(56')\ \ [b,b'\circ b'']+[b',b''\circ b]+[b'',b\circ b']=0,\label{eq:5.14}
\end{eqnarray} 
for $a,a',a''\in A$ and $b,b'\in B$.  \medskip

First we establish a general lemma
for any $\hbox{\rm C}_2$-graded Lie algebra.  \bigskip
  
\begin{lem}\label{Lem:5.15} Suppose that $\fa = A \oplus B$ is the coordinate algebra
of a $\hbox{\rm C}_2$-graded Lie algebra.    For $b,b' \in B$,  suppose that there exist
elements $a_1,a_2,a_3,a_4\in A$ such that
$b=\hf [a_1,a_2]$ and $b'=\hf [a_3,a_4]$.    Then 
$D_{b,b'}=[D_{a_1,a_2},D_{a_3,a_4}]$.
Hence,  $D_{b,b'}$ restricted to the Jordan algebra $(A, \cdot)$ is an inner derivation.
\end{lem}

\pf  We know that $D_{b,b'}$ is an inner derivation of the Jordan algebra $(\fa^+, \cdot)$.
What this result asserts is that $D_{b,b'}$ acts as an inner derivation 
of the Jordan algebra $(A, \cdot)$.    

By  \eqref{eq:5.3}, we have
$[b,a]=\hf[[a_1,a_2],a]]= 2D_{a_1,a_2}a$ and 
$[b',a]=2D_{a_3,a_4}a$. Then, by the same reason,
$[b,D_{a_3,a_4}a]=2D_{a_1,a_2}D_{a_3,a_4}a$
and
$[b',D_{a_1,a_2}a] =2D_{a_3,a_4}D_{a_1,a_2}a$. Therefore,
by \eqref{eq:5.11},

\begin{eqnarray*} D_{b,b'}a&=&\frac{1}{4}\big ([b,[b',a]]-[b',[b,a]]\big) \\ 
&=&\frac{1}{2}\big([b,D_{a_3,a_4}a]-[b',D_{a_1,a_2}a]\big)
\\ &=&D_{a_1,a_2}D_{a_3,a_4}a -D_{a_3,a_4}D_{a_1,a_2}a
\\ &=&[D_{a_1,a_2},D_{a_3,a_4}]a. \qed
\end{eqnarray*} \medskip
 
Next we impose the assumptions that the Lie algebra is  
$(\hbox{\rm C}_2,G)$-graded  and satisfies the division property.  
  \bigskip
  
\begin{lem} \label{Lem:5.16} 
Let  $\fa = A\oplus B$ be the coordinate algebra of a division $(\hbox{\rm C}_2,G)$-graded Lie algebra.
Assume $0\neq a,a'\in A$ and $0\neq b\in B$ are homogeneous. Then 
\begin{itemize}
\item[{\rm (i)}]  $a\circ a'\neq 0$ or $[a,a']\neq 0$;
\item[{\rm (ii)}]  $a\circ b\neq 0$ or $[a,b]\neq 0$.
\end{itemize}
\end{lem}

\pf    For (i), suppose that $a\circ a'=0=[a,a']$. Then by \eqref{eq:5.1},  we have
$$ a\circ (a'\circ a'') =[a,[a'',a']].
$$ But then substituting  ${a'}^{-1}$ for $a''$, gives 
$4a=0$, a contradiction.

For (ii), suppose that $a\circ b=0=[a,b]$. Then by \eqref{eq:5.9},   we have
$$ [a,[b,a']]=-b\circ (a\circ a').
$$ Letting $a' = a^{-1}$ gives
$[a,[b,a^{-1}]]=-4b$. But, by \eqref{eq:5.8},  we have
$$ 0=4D_{a,a^{-1}}b =[a,[a^{-1},b]]-[a^{-1}[a,b]].
$$ Hence,
$-4b=[a,[b,a^{-1}]]=[a^{-1},[a,b]]=0$,  a contradiction.
\qed  
\bigskip

Recall that for a $\hbox{\rm C}_2$-graded Lie algebra, the skew product  on $B$ 
may be arbitrarily defined.    Here we will make a precise definition of that
skew product when the Lie algebra is a Lie $G$-torus of type C$_2$ .  This then will enable us 
to determine the structure of the corresponding coordinate algebra 
$\fa=A\oplus B =\bigoplus_{g\in B}\ (A^g\oplus B^g)$.     The Lie $G$-torus condition  
forces  $\dim_\F\fa^g\leq 1$
for all $g\in G$, where $\fa^g=A^g\oplus B^g$,
and so this implies the helpful fact that $A^g=0$ or $B^g=0$. \bigskip

\begin{lem}\label{Lem:5.17} 
Let  $\fa = A\oplus B$ be the coordinate algebra of a Lie $G$-torus of type $\hbox{\rm C}_2$.
Set $B_0:=\{b\in B\mid [b,A]=0\}$.   Then, 
$$B=[A,A]\oplus B_0.$$ Moreover,
$[A,A]$ and $B_0$ are graded, and for any homogeneous element $b\in [A,A]$, there exist  homogeneous elements
$a_1,a_2\in A$ such that
$b=[a_1,a_2]$.
\end{lem}

\pf   Clearly $[A,A]$ and $B_0$ are graded spaces. 
Suppose that $0\neq  b\in {B^g} \setminus B_0$. Then, there exists some $a\in A^h$ such that
$[b,a]\neq 0$. Note that $a^{-1}\circ [b,a]\in A^g=0$
by one-dimensionality. Hence, by Lemma \ref{Lem:5.16}, we have $[a^{-1}, [b,a]]\neq
0$.  Since $B^g$ is one-dimensional, there must exist  some $\vartheta \in \F$ such that
$b=\vartheta[a^{-1}, [b,a]]$.   Thus if one takes $a_1=\vartheta a^{-1}$ and $a_2=[b,a]$, 
then $b=[a_1,a_2]$, and it follows that $B=[A,A]+B_0$.

Suppose that $b=[a_h,a_k]\in B_0$,  where 
$0\neq a_h\in A^h$ and $0\neq a_k\in A^k$. Then $[b,a_h^{-1}]=0$, and hence $b\circ a_h^{-1}\neq 0$ by
Lemma \ref{Lem:5.16}.   Thus, by the one-dimensionality of $A^k$, there exists some  $\xi \in F$ such that
$a_k=\xi b\circ a_h^{-1}$.   We apply identity \eqref{eq:5.7} with $a=a_h$ and $a'=\xi a_h^{-1}$ to obtain  $0=b\circ[a_h,\xi a_h^{-1}] +[b\circ \xi a_h^{-1},a_h]=[b\circ \xi a_h^{-1},a_h]
=[a_k,a_h]=b$, 
and so $[A,A]\cap B_0=0$.
\qed

\medskip

Set $S:=\supp\fa$, $S_+:=\supp A$,  and $S_-:=\supp B$.   Then 
$S=S_+\sqcup S_-$, which is a disjoint union by the one-dimensionality condition on
the graded spaces of $\fa$.
\bigskip

\begin{lem}\label{Lem:5.18} 
Let  $\fa= A\oplus B$ be the coordinate algebra of a Lie $G$-torus of type $\hbox{\rm C}_2$.  Then
the following are equivalent:
\begin{itemize}
\item[{\rm (i)}]  $S_+$ is a subgroup of $G$;
\item[{\rm (ii)}]  $[A,B]=0$; 
\item[{\rm (iii)}] $[A,A]=0$;   
\item[{\rm (iv)}]  The product $\cdot$  on $A$ coincides with the product from $\fa$;
\item[{\rm (v)}] 
$(A,\cdot)$ is associative.
\end{itemize}
\end{lem}

\pf
It follows from  Lemma \ref{Lem:5.17} that $[A,B]=0$ if and only if $[A,A]=0$.
In this case, the product of $\fa$ and $\cdot$
coincide on $A$, and by \eqref{eq:5.3},
$(A,\cdot)$ is associative. 
Then, by the division property,
the support $S_+$ of the division graded associative algebra $(A,\cdot)$ is a subgroup of $G$.
Thus,  we only need to prove that $[A,B]=0$ if $S_+$ is a subgroup.
For $g\in S_+$,
$h\in S_-$,  and $b\in B^h$, we have
$[A^g,b]\subseteq A^{g+h}$. But since $S_+$ is a subgroup, $g+h\notin S_+$. Hence, $A^{g+h}=0$, and so
$[A,B]=0$.
\qed 
\bigskip

At this juncture, it is convenient to divide our considerations into two cases, namely, 
\begin{itemize}
\item[{\rm (i)}]  $S_+$ is a subgroup;
\item[{\rm (ii)}] $S_+$ is not a subgroup.
\end{itemize} \m 

Case (i): When $S_+$ is a subgroup, we will show that  if we set $[B,B]=0$, 
then $\fa=\fa^+$ is a Clifford $G$-torus. 
By Lemma \ref{Lem:5.18}, we know that 
the product of $\fa$ coincides with $\cdot$
on $A$,  and 
$(A,\cdot)$ is a commutative, associative algebra.
Now by \eqref{eq:5.9}, we have
$a'\cdot(a\cdot b)\cdot a'=(a\cdot a')\cdot b$, i.e.,
$B$ is a graded $(A,\cdot)$-module by the action $\cdot$. Also,  \eqref{eq:5.11}  and \eqref{eq:5.2} imply that
$$ b\circ (b'\circ a) =b'\circ(b\circ a) =a\circ(b\circ b'),
$$ and so
$\circ$ (and also $\cdot$) defines a symmetric $A$-bilinear form on $B$.
Note that the form is nondegenerate  by the division property.
Thus, if we specify that $[B,B]=0$, then  $\fa=\fa^+$ is a Jordan algebra 
of a symmetric bilinear form $\cdot$
on $B$ over $A$. (Observe we did not use the fact that
$\fa^+$ is a Jordan algebra in showing this.)     Therefore,  we have that 
$\fa$ is a Clifford $G$-torus.   \medskip

Case (ii): When $S_+$ is not a subgroup, we define a linear map
$\vvp: [A,A]\longrightarrow \inder (A,\circ)$  into the inner derivations of the Jordan algebra
$(A,\circ)$, by requiring that $\vvp(b) \in  \inder (A,\circ)$ be given by
$\vvp(b)(a)=[b,a]$ for all $b\in [A,A]$ and $a\in A$.  By  \eqref{eq:5.3}, the image of $\vvp$ is indeed in $\inder
(A,\circ)$,  and $\vvp$ is surjective. Moreover, by Lemma \ref{Lem:5.17}, $\vvp$ is injective. 
We will use  the bijection $\vvp$ to define a
skew product $[b,b']$ for $b,b'\in B$ in the following way.    If $b,b'\in [A,A]$, then by 
Lemma \ref{Lem:5.15},
$D_{b,b'}\in\inder (A,\circ)$.   Hence, there is a unique element in $[A,A]$, denote it $[b,b']$, 
so that 
$[b,b']:=\vvp^{-1}(4D_{b,b'})$.    If $b\in B_0$ or $b'\in B_0$, we specify that  $[b,b']:=0$.   Thus we have the
following relation:
\begin{equation}\label{eq:5.19} 4D_{b,b'}a=[[b,b'],a]
 \end{equation}
(see \eqref{eq:5.11}  for $b\in B_0$ or $b'\in B_0$), which is a well-known identity for an associative algebra.

Now we can prove that $\fa$ is associative in exactly the same way as in \cite[pp.~105-111]{Se}.  Indeed, we have established all the  properties needed in Seligman's argument  to show that  the associative law holds in $\fa$ except for the
simplicity of our Lie algebra $\scl$.    However, 
a centerless Lie $G$-torus is
graded simple; that is, it has no nontrivial graded ideals (see \cite [Lem.~4.4]{Y3}), and simplicity may be replaced by graded simplicity
with no harm to the argument.
For the convenience of the reader, in the paragraphs to follow we present an alternative proof
of associativity  which differs somewhat from Seligman's original argument.      
The ultimate conclusion of Case (ii) will be that $\fa$ is an associative $G$-torus with
graded involution.

The $D$-mappings are derivations   
not only relative to  the symmetric product $\circ$  but also
relative to the skew product $[\cdot, \cdot]$ we have defined above.
First we prove the following claim  about $D_{B,B}B$.
(Seligman proved this claim using the Lie algebra $\scl$,
but we prove it using just the coordinate algebra $\fa$.)

\begin{claim}\label{claim:5.20}
When $S_+ = \supp A$ is not a subgroup, then 
$D_{B_0,B_0}B_0=0=D_{B_0,B_0}B$,
and therefore the following hold:
\begin{itemize} 
\item[{\rm (i)}] $D_{B_0,B_0}=0=D_{B_0,B},$
\item[{\rm (ii)}] $D_{B,B}B_0=0,$ 
\item[{\rm (iii)}] $D_{B,B}=D_{[A,A],[A,A]}\subset D_{A,A}.
$  \end{itemize}
Thus 
$D_{B,B}B\subseteq  D_{A,A}B\subseteq D_{A,A}[A,A]\subseteq [A,A]$.
\end{claim}

\pf  {F}rom Lemma \ref{Lem:5.17} 
we have
$$
D_{B,B}=D_{[A,A]\oplus B_0,[A,A]\oplus B_0}.
$$
By \eqref{eq:5.4}, 
$D_{[A,A],B_0}=0$, so the above becomes 
$$
D_{[A,A],[A,A]}+D_{B_0,B_0}\subseteq D_{A,A}+D_{B_0,B_0},
$$
again using  \eqref{eq:5.4}.
Now  \eqref{eq:5.11} implies  $D_{B_0,B_0}A =0$, and so $D_{B_0,B_0}[A,A]$ $=0$.
Hence,
$D_{B_0,B_0}B=D_{B_0,B_0}B_0\subseteq B_0$.
Since $D_{B,B_0}=D_{B_0,B_0}$ by the above, we have
 $D_{B,B_0}B=D_{B_0,B_0}B_0$; however, 
 \eqref{eq:5.8} gives
$D_{A,A}B_0=0$, so that $D_{B,B}B_0=D_{B_0,B_0}B_0$ as well. 

Now set
$$
K:=B\circ D_{B_0,B_0}B_0
+ D_{B_0,B_0}B_0.
$$
We claim $K$ is an ideal of  $\fa^+$.
To see this, note that 

\begin{eqnarray*}
B\circ D_{B_0,B_0}B_0
&\subseteq& D_{B_0,B_0}(B\circ B_0) +(D_{B_0,B_0}B_0)\circ B_0\\
&\subseteq& D_{B_0,B_0}A+(D_{B_0,B_0}B_0)\circ B_0
\subseteq (D_{B_0,B_0}B_0)\circ B_0,
\end{eqnarray*}

\noindent and so
$$
B\circ(B\circ D_{B_0,B_0}B_0)\subseteq B\circ(B_0\circ D_{B_0,B_0}B_0),
$$
which 
is contained in
$$
B_0\circ (B\circ D_{B_0,B_0}B_0)+D_{B,B_0}D_{B_0,B_0}B_0
\subseteq
B_0\circ (B_0\circ D_{B_0,B_0}B_0)+D_{B_0,B_0}B_0. 
$$
Now $B_0\circ (B_0\circ D_{B_0,B_0}B_0)+D_{B_0,B_0}B_0$  lies in 
$(D_{B_0,B_0}B_0)\circ (B_0\circ B_0)+D_{B_0,B_0}B_0$, 
since 
\begin{eqnarray*} B_0\circ (B_0\circ D_{B_0,B_0}B_0)
&=&D_{D_{B_0,B_0}B_0, B_0}B_0+(D_{B_0,B_0}B_0)\circ (B_0\circ B_0) \\
&& \qquad \subseteq D_{B_0,B_0}B_0+(D_{B_0,B_0}B_0)\circ (B_0\circ B_0).
\end{eqnarray*}
But
\begin{eqnarray}\label{eq:5.21}
(D_{B_0,B_0}B_0)\circ A&\subseteq& D_{B_0,B_0}(B_0\circ A)+B_0\circ D_{B_0,B_0}A=D_{B_0,B_0}(B_0\circ A) \nonumber \\
&\subseteq&  D_{B_0,B_0}B=D_{B_0,B_0}B_0, \end{eqnarray}
so the above shows
$$
B\circ (B\circ D_{B_0,B_0}B_0)\subset D_{B_0,B_0}B_0.
$$
Thus to verify that $K$ is an ideal, it suffices to show 
$A\circ (B\circ D_{B_0,B_0}B_0)\subseteq  B\circ D_{B_0,B_0}B_0$.
But 
$B\circ D_{B_0,B_0}B_0=B_0\circ D_{B_0,B_0}B_0$,  so
the required inclusion follows from \eqref{eq:5.2}  and \eqref{eq:5.21}.  
Consequently,  $K$ is an ideal in $\fa^+$ as claimed.

Clearly $K$ is graded, and $\fa^+$ is graded simple.
Hence, $K=0$ or $K=\fa^+$.
Suppose $K=\fa^+$.  Then $A=B\circ D_{B_0,B_0}B_0$.
But 
we have
$[B,B\circ D_{B_0,B_0}B_0]=[[A,A],B\circ D_{B_0,B_0}B_0]
=[[A,A],(D_{B_0,B_0}B_0)\circ B_0]$, and
by \eqref{eq:5.3}, this is contained in
$$
D_{A,A}((D_{B_0,B_0}B_0)\circ B_0)=0,
$$
since $D_{A,A}B_0=0$.
That is, we have
$[B,A]=[B, B\circ D_{B_0,B_0}B_0]=0$,
which is not our case.
Therefore $K=0$, and all the statements 
in Claim \ref{claim:5.20} follow.    
\qed
 
 \bigskip
 Now, we want to show  that the associator
$(\alpha, \alpha',\alpha''):=(\alpha\alpha')\alpha''-\alpha(\alpha'\alpha'')=0$,
and for this purpose, it suffices to verify that the following two identities
\begin{eqnarray}
&&(\alpha\circ\alpha')\circ\alpha''-\alpha\circ(\alpha'\circ\alpha'')
=[\alpha,[\alpha',\alpha'']]-[[\alpha,\alpha'],\alpha'']\qquad  \label{5.22} \\
&&[\alpha,\alpha']\circ\alpha''-\alpha\circ[\alpha',\alpha'']
=[\alpha,\alpha'\circ\alpha'']-[\alpha\circ\alpha',\alpha'']. \qquad   \label{5.23}
\end{eqnarray}
\noindent are satisfied for all $\alpha,\alpha',\alpha''\in\fa$.
By \eqref{eq:5.1}, relation (\ref{5.22})  holds for $\alpha,\alpha',\alpha''\in A$.
Note that from the definition of $D_{\alpha,\alpha'}$ in \eqref{eq:der}, we have
$D_{\alpha'',\alpha}\alpha'
=\alpha''\circ(\alpha\circ\alpha')-\alpha\circ(\alpha''\circ\alpha')$,
which is the left-hand side of  (\ref{5.22}).  In addition, we have (\ref{5.22})  for
$\alpha,\alpha''\in B$ and $\alpha'\in A$ by \eqref{eq:5.11}.
Interchanging $a$ and $a'$ in \eqref{eq:5.9}
and subtracting the two relations
gives 
(\ref{5.22})  for
$\alpha,\alpha''\in A$ and $\alpha'\in B$.
Thus the cases remaining to establish (\ref{5.22})  are:

\begin{itemize}
\item[{\rm (a)}] \quad $(a\circ b)\circ b'-a\circ(b\circ b')
= [a,[b,b']]-[[a,b],b']$
\item[{\rm (b)}] \quad $(a\circ a')\circ b-a\circ(a'\circ b)
=[a,[a',b]]-[[a,a'],b]$ 
\item[{\rm(c)}] \quad $(b\circ b')\circ b''-b\circ(b'\circ b'')
=[b,[b',b'']]-[[b,b'],b'']$
\end{itemize}
for all $a,a',a''\in A$ and $b,b',b''\in B$.

For (a), starting with \eqref{eq:5.2}, we have
\begin{eqnarray*}
2a\circ(b\circ b')
&=& b\circ(b'\circ a)+b'\circ(b\circ a) +[b,[b',a]]+[b',[b,a]]\\
&=&  b\circ(b'\circ a)+b'\circ(b\circ a) \\ && \hspace{.6 truein}
+[[b,b'],a]+2[b',[b,a]]
\quad\text{by \eqref{eq:5.19} and \eqref{eq:5.11}}\\
&=& 2 b\circ(b'\circ a)+2[b',[b,a]] \\
&& \hspace{.6 truein} \text{by \eqref{eq:5.19} and the definition of $D_{b,b'}$.}
\end{eqnarray*}
Hence,  $a\circ(b\circ b')= b\circ(b'\circ a)+[b',[b,a]]$, and  thus, 
\begin{eqnarray*}(a\circ b)\circ b'-a\circ(b\circ b') &=&(a\circ b)\circ b'-b\circ(b'\circ a)-[b',[b,a]]\\
&=&-[[b,b'],a]-[b',[b,a]]  
\quad\text{by \eqref{eq:5.19},} \\
&=& [a,[b,b']]-[[a,b],b'].
\end{eqnarray*}
 
Now applying \eqref{eq:5.9},  we see  that  (b) holds once  we check that
\begin{equation} \label{eq:5.24}
[[a,a'],b]=[a,[a',b]]-[a',[a,b]] . \end{equation}
 Recall that $[[a,a'],b]$ is a unique element in $B$ so that
$[[[a,a'],b],a'']=4D_{[[a,a'],b]}a''$ for all $a''\in A$.
But then 
\begin{eqnarray*} 
[[[a,[a',b]]-[a',[a,b]],b],a'']
&=&[[[a,[a',b]],b],a'']-[[[a',[a,b]],b],a'']\\
&=&4D_{[[a,[a',b]],b]}a''-4D_{[a',[a,b]],b]}a''
\,\text{by (\ref{5.23})}\\
&=&4D_{[[a,a'],b]}a''
\quad\text{by \eqref{eq:5.4}}.
\end{eqnarray*} As a result,  we obtain \eqref{eq:5.24}.  
 
 \m  Finally for (c),  we observe that 
 $(b\circ b')\circ b''-b\circ(b'\circ b'')=D_{b'',b}b'$
is in $[A,A]$  by Claim \ref{claim:5.20}. 
So, by the decomposition of $B$ in Lemma \ref{Lem:5.17},
it suffices to show that, for all $a\in A$,
$[D_{b'',b}b',a]+[[[b,b'],b''],a]=[[b,[b',b'']],a]$,
or that
$D_{b'',b}[b',a]-[b'D_{b'',b}a]
+D_{[b,b'],b''}a=D_{b,[b',b'']}a$.  
By \eqref{eq:5.11}, this amounts to verifying that 
\begin{eqnarray*}
&&\hspace{-.5 truein} D_{b'',b}[b',a]-[b',D_{b'',b}a]
+[[b,b'],[b'',a]]
-[b'',[[b,b'],a]]\\
&& \hspace {.8truein} =[b,[[b',b''],a]]
-[[b',b''], [b,a]], \\
\end{eqnarray*}  or  
$$
D_{b'',b}[b',a]-[b',D_{b'',b}a]
+D_{b,b'}[b'',a]
-[b'',D_{b,b'}a]=[b,D_{b',b''}a]
-D_{b',b''}[b,a],
$$  or 
$$
[D_{b'',b}b'+D_{b,b'}b''+D_{b',b''}b,a]=0.
$$
But 
$D_{b'',b}b'+D_{b,b'}b''+D_{b',b''}b=0$
(simply by the definition of the inner derivation), and hence (c) is established.

Now for identity (\ref{5.23}),
interchanging $\alpha$ and $\alpha''$ in \eqref{eq:5.2}
will show that this identity holds for $\alpha,\alpha',\alpha''\in A$.
Also,  \eqref{eq:5.7}  implies  the case with 
$\alpha,\alpha'\in A$ and  $\alpha''\in B$,
and  the case with
$\alpha',\alpha''\in B$ and $\alpha\in A$.
What remains to be checked is that the following
equations hold for all $a,a',a'' \in A$ and $b,b',b'' \in B$:
\smallskip
\begin{itemize}
\item[{\rm (d)}]  \quad ${[a,b]}\circ a'-a\circ[b,a']
=[a,b\circ a']-[a\circ b,a']$
\item[{\rm (e)}]  \quad
${[b,a]}\circ b'-b\circ[a,b']
=[b,a\circ b']-[b\circ a,b']$
\item[{\rm (f)}] \quad ${[b,b']}\circ a-b\circ[b',a]
=[b,b'\circ a]-[b\circ b',a]$
\item[{\rm (g)}] \quad ${[b,b']}\circ b''-b\circ[b',b'']
=[b,b'\circ b'']-[b\circ b',b'']$. \end{itemize}
\smallskip

Reversing the roles of $a$ and $a'$ in \eqref{eq:5.7} and adding gives
$$
[b\circ a, a']+[b\circ a',a]+a\circ[b,a']+a'\circ[b,a]=2[b,a\circ a'].
$$
Then by \eqref{eq:5.5}, we have (d).

In view of \eqref{eq:5.13},
verification of (e) reduces to showing 
$$
[b\circ a,b']=[a,b\circ b']+[b,b'\circ a],
$$
or, since all terms belong to $[A,A]$, that for all $a''\in A$,
$$
[[b\circ a, b'],a'']+[[b'\circ a,b],a'']
=[[a,b\circ b'],a''].
$$
By definition, the left side equals $4D_{b\circ a, b'}a''+4D_{b'\circ a, b}a''$,
which by \eqref{eq:5.10} is
$4D_{a,b\circ b'}a''$. 
The relation in (e) now follows from \eqref{eq:5.3}.

To establish (f),
it will suffice by the decomposition $B = [A,A] \oplus B_0$
to treat two separate cases:

(1) $b'=b_0\in B_0$;

(2) $b'=[a',a'']$ for some $a',a''\in A$.

Now when  $b'=b_0$, we have
$[b,b']=0$,
and our relation reduces to showing
$$
[b\circ b_0, a]=[b, b_0\circ a].
$$
Since  both members are in $[A,A]$,
it suffices to prove that
$$
[[b\circ b_0,a],a']=[[b,b\circ a],a'].
$$ for $a'\in A$.  By \eqref{eq:5.3}, the left side is $D_{b\circ b_0,a}a'$,
while the right equals $D_{b, b_0\circ a}a'$ by 
the definition of $[\cdot,\cdot]$ on $B$.
Equation \eqref{eq:5.10} shows that the difference of the two is
$D_{b\circ a, b_0}a'=0$, 
since $D_{B,B_0}=0$ by Claim \ref{claim:5.20}.
Thus (f) holds if $b'\in B_0$.

Suppose then that  $b'=[a',a'']$, where $a',a''\in A$.
Here (f) reads:

\begin{equation}\label{eq:5.25}
 \hspace{.15 truein} [b,[a',a'']]\circ a-b\circ [[a',a''],a]
=[b,[a',a'']\circ a]-[b\circ[a',a''],a]. \end{equation}

\noindent We have
$[b,[a',a'']]=(b\circ a')\circ a''+[[b,a'],a'']-b\circ (a'\circ a'')$ from (\ref{5.22}).
Substitution shows the first term on the left-hand side of \eqref{eq:5.25}
to be:
$$
((b\circ a')\circ a'')\circ a+[[b,a'],a'']\circ a-(b\circ (a'\circ a''))\circ a.
$$
By \eqref{eq:5.3}, the second term on the left of \eqref{eq:5.25} is
$$
-b\circ D_{a',a''}a=-D_{a',a''}(b\circ  a)+(D_{a',a''}b)\circ a
$$
by the derivation property.  This expression is equal to
$$
-[a',[a'',b\circ a]]+[a'',[a',b\circ a]]+[a',[a'',b]]\circ a-[a'',[a',b]]\circ a
$$
by \eqref{eq:5.8}. 
Hence, the left-hand side of  \eqref{eq:5.25}
becomes 
\begin{equation}\label{eq:5.26}
-[a',[a'',b\circ a]]+[a'',[a',b\circ a]]
+((b\circ a')\circ a''-b\circ (a'\circ a'')+[a',[a'',b]])\circ a,
\end{equation}
which is equal to
$-[a',[a'',b\circ a]]+[a'',[a',b\circ a]]$  by \eqref{eq:5.9}.
Thus (f) reduces to showing
\begin{equation} \label{eq:5.27} 
 \ \
-[a',[a'',b\circ a]]+[a'',[a',b\circ a]]
=[b,[a',a'']\circ a]-[b\circ[a',a''],a].
\end{equation}
 
All terms in \eqref{eq:5.27}  lie in $[A,A]$,
so it is sufficient to show, for all $a'''\in A$, that 
$$
[-[a',[a'',b\circ a]]+[a'',[a',b\circ a]], a''']=[[b,[a',a'']\circ a]-[b\circ[a',a''],a], a'''].
$$
By \eqref{eq:5.3} and our definition of $[\cdot,\cdot]$ on $B$,
this identity is equivalent to
$-D_{a',[a'',b\circ a]}+D_{a'',[a',b\circ a]}
=D_{b,[a',a'']\circ a}-D_{b\circ[a',a''],a}$.
Then \eqref{eq:5.10} can be quoted to give  
$-D_{a',[a'',b\circ a]}+D_{a'',[a',b\circ a]}
=D_{a\circ b,[a',a'']}$,
and then (f)  is a direct consequence of  \eqref{eq:5.4}.

Finally,  we tackle (g).
As in  case  (1),   where $b'\in B_0$, we must show
$[b,b'\circ b'']+[b'',b\circ  b']=0$
for all $b,b''\in B$.
However,  this is immediate from \eqref{eq:5.14}.
For  case (2), where $b'=[a,a']$, 
we use (b) of (\ref{5.22})  to substitute for
$[[a,a'],b]$ and $[[a,a'],b'']$, then apply \eqref{eq:5.9} to obtain
$[b,b']=(a\circ b)\circ a'-a\circ (b\circ a')$
and
$[b,b'']=a\circ (b''\circ a')-(a\circ b'')\circ a'$.
Now showing (g) reduces, by (\eqref{eq:5.14}, to showing
$$
[b',b''\circ  b]+[b,b']\circ b''-[b',b'']\circ b=0
$$
for $b'=[a,a']$, or that
\begin{equation}\label{eq:5.28}
  [[a,a'],b''\circ  b]+[b,[a,a']]\circ b''-[[a,a'],b'']\circ b=0.
\end{equation}
Now $[[a,a'],b''\circ b]=D_{a,a'}(b''\circ b)$ by \eqref{eq:5.3},
and this says
\begin{eqnarray*} 
&&\hspace{-.2 truein}(D_{a,a'}b'')\circ b+(D_{a,a'}b)\circ b'' \\
&& \hspace{.2 truein} =[a',[b'',a]]\circ  b-[a,[b'',a']]\circ b
+[a',[b,a]]\circ b''-[a, [b,a']]\circ b''
\end{eqnarray*}
by \eqref{eq:5.8}
Thus, the left side of \eqref{eq:5.28} becomes   
\begin{eqnarray*}
&&\hspace{-.3 truein} \bigg([a',[b'',a]]\circ  b-[a,[b'',a']]-a\circ (b''\circ a')+(a\circ b'')\circ a'\bigg)\circ b\\
&&+
\bigg([a',[b,a]]\circ b''-[a, [b,a']]+(a\circ b)\circ a'-a\circ (b\circ a')\bigg)\circ b'',
\end{eqnarray*}
which is 0 by (\ref{5.22}).  Thus,  (g) is proved, and  
we have the desired conclusion that $\fa$ is associative. 
\bigskip 

\begin{rem}\label{Rem:5.23}
{\rm When $S_+ = \supp A$ is not a subgroup of $G$,
 we have $[A,B_0]=0=[B,B_0]$, so
$B_0$  is contained in the center of the associative algebra 
$\fa=A\oplus B$.}  \end{rem} \medskip

\begin{rem}\label{Rem:5.24}
{\rm   When $S_+$ is not a subgroup,
then we claim that  $G=\langle S_+\rangle$ is forced.
To see this, let  $G':=\langle S_+\rangle$ and
suppose that $G'\neq G$.
Then, there exists some  $g\in S_-\setminus G'$.
Take  $0\neq b\in\fa^g=B^g$.
Then, $[A,b]=0$.
Indeed, let $0\neq a\in\fa^{h}=A^{h}$.
Then, $[a,b]\in A^{h+g}$,
but $h+g\notin G'$ since $h \in G'$.
Hence, by Lemma \ref{Lem:5.16}, $b\circ a\neq 0$ for any homogeneous $0\neq a\in A^h$.
For the same reason,
$(b\circ a)\circ a'\neq 0$ for any homogeneous  $0\neq a'\in A$
since $b\circ a\in B^{h+g}$ and $h+g\notin G'$.
However, there exist nonzero homogeneous elements $a,a'\in A$
such that $a\circ a'=0$ since $S_+$ is not a subgroup.
This is a contradiction
by \eqref{eq:5.9};
$$
(b\circ a)\circ a'=b\circ (a\circ a')
$$
for any $a,a'\in A$.

Consequently,  for an associative $G$-torus $\fa=A\oplus B$ 
with involution, if  the set  $S_+$ does not generate $G$,
then $S_+$ is a subgroup,  and
$\fa$ becomes a Clifford $G$-torus upon defining  $[B,B]=0$.}
\end{rem}
\medskip

Combining all our results,  we arrive at our main theorem. 
\bigskip

\begin{thm}\label{Thm: 5.13}
A centerless Lie $G$-torus $\scl$ of type $\hbox{\rm C}_r$ is isomorphic
to  a Lie algebra
 $\mathfrak{sp}_{2r} (\fa)$,
where $\fa$ is:
\begin{itemize}
\item  an associative $G$-torus
with graded involution if $r\geq 4$,
\item an alternative $G$-torus 
with graded involution whose symmetric elements are in the nucleus of 
$\fa$  if $r=3$,
\item an associative $G$-torus
with graded involution or a Clifford $G$-torus if $r=2$.
\end{itemize}
\end{thm}

This result generalizes the classification of the core of extended affine Lie algebras
of type $\hbox{\rm C}_r$ in \cite{AG},  as  the core is a Lie torus, i.e., 
a Lie $G$-torus for $G=\mathbb Z^n$.
In this case one can use more concrete terminology in describing $\scl$.  \bigskip

\begin{cor}\label{cor:5.14} {\rm (Compare \cite[Thm.~4.87]{AG}.)} 
A centerless Lie torus $\scl$ of type $\hbox{\rm C}_r$ is isomorphic to a Lie
algebra  $\mathfrak{sp}_{2r} (\fa)$,
where $\fa$ is:
\begin{itemize}
\item  a quantum torus with graded involution if $r\geq 4$,
\item a quantum torus with graded involution or 
an octonion torus with standard involution if $r=3$,
\item a quantum torus with graded involution or a Clifford torus if $r=2$.
\end{itemize}
\end{cor}

\medskip

\noindent  {\small \sc Department of Mathematics, University of 
Wisconsin, Madison, WI 53706-1388 USA

\noindent  {\tt benkart@math.wisc.edu}

\medskip
 
\noindent  {\small \sc Department of Mathematics, North Dakota State University,  Fargo, ND 58105-5075  USA

\noindent  {\tt yoji.yoshii@ndsu.nodak.edu}


\smallskip
 
\begin{thebibliography}{888}


 \bibitem[AABGP]{AABGP} B.N.~Allison, S.~Azam, S.~Berman, Y.~Gao,
 and A.~Pianzola, {\it Extended Affine Lie Algebras and Their
 Root Systems}, Mem. Amer. Math. Soc. \textbf {126}  vol.~603,
 Providence, R.I.  1997.

 \bibitem[ABG1]{ABG1}  B.N.~Allison, G.~Benkart, and Y.~Gao, Central
extensions of Lie algebras graded by finite root systems, Math. Ann. 
\textbf{ 316} ( 2000), 499-527.


\bibitem[ABG2]{ABG2}  B.N.~Allison, G.~Benkart, and Y.~Gao,
{\it Lie Algebras Graded by the Root Systems BC$_r$, $r \geq 2$}, 
Mem. Amer.  Math.  Soc.  \textbf {158} vol.~751, Providence, R.I. 
2002.

 
\bibitem[ABG3]{ABG3}  B.N.~Allison, G.~Benkart, and Y.~Gao,
Extended affine Lie algebras of type $\hbox{\rm BC}_r$, $r \geq 3$,
preprint.  


\bibitem[AFY]{AFY}  B.~Allison, J.~Faulkner,  and Y. Yoshii, 
Structurable tori, to appear.   

\bibitem[AG]{AG}  B.N.~Allison and Y.~Gao,
The root system and the core of an extended affine Lie algebra, 
Selecta Math. (N.S.)
\textbf{7} 
(2001),  149--212. 

\bibitem[AY]{AY}  B.N.~Allison and Y. Yoshii, 
Structurable tori and extended affine Lie algebras of type 
$\hbox{\rm BC}_1$, 
Pure Appl. Algebra, no. 2-3,
\textbf{184} (2003)
105--138. 

\bibitem[BN]{BN} G.~Benkart and E.~Neher, The centroid of
extended affine and root graded Lie algebras, J. Pure and
Appl. Alg., to appear.


\bibitem[BS]{BS} G.~Benkart and O.~Smirnov, Lie algebras 
graded by the root system BC$_1$,   J.  Lie Theory \textbf{13} 
(2003), 91--132.

\bibitem[BZ]{BZ}
G.~Benkart and E.~Zelmanov, Lie algebras graded by finite root 
systems and intersection matrix algebras, Invent. Math. \textbf{126} 
(1996), 1--45. 

\bibitem[BGK]{BGK} S.~Berman, Y.~Gao, Y.~Krylyuk,  
Quantum tori and the structure of elliptic quasi-simple Lie algebras,
J. Funct. Anal. 
\textbf{135}
(1996)  339--389.

 \bibitem[BGKN]{BGKN}  S.~Berman, Y.~Gao,  Y.~Krylyuk, and
E. Neher, The alternative torus and the structure of elliptic 
quasi-simple Lie algebras of type $A_2$, Trans. Amer. Math. Soc. 
\textbf{347} (1995),  4315--4363.

\bibitem[BM]{BM} S.~Berman and R.V.~Moody,  Lie algebras
graded by finite root systems and the intersection matrix algebras 
of Slodowy, Invent. Math.  \textbf {108} (1992), 323--347.

\bibitem[B]{B} N.~Bourbaki,  {\it Groupes et Alg\`ebres de Lie,
Chapitres 4--6}, Masson, Paris 1981.  

\bibitem[F]{F} J.R.~Faulkner, Lie tori of type BC$_2$ and structurable quasitori,
preprint. 
 
\bibitem[J]{J}  N. Jacobson, 
{\it Structure and Representations of Jordan Algebras},
Amer. Math. Soc. Colloq. Publ. vol.\textbf{~39},
Providence R.I. 
1968.

\bibitem[M]{M} K. McCrimmon,
{\it A Taste of Jordan Algebras}, 
Universitext, Springer-Verlag
New York 2004.

 
\bibitem[MN]{MN}  K.~McCrimmon and E.~Neher,  
Coordinatization of triangulated Jordan systems,
 J.~Algebra
\textbf{114}
(1988)  411--451.

\bibitem[N1]{N1}  E.~Neher, Lie tori,  C.~R.~Math.~Acad.~Sci.~Soc.~R.~Can.  \textbf{26}  (2004),  no.~3, 84--89.
\bibitem[N2]{N2}  E.~Neher, Extended affine Lie algebras, C.~R.~Math.~Acad.~Sci.~Soc.~R.~Can. \textbf{26} (2004), no.~3,  90--96.

\bibitem[N3]{N3}  E.~Neher, 
Quadratic Jordan superpairs covered by grids.  J. Algebra  \textbf{269}  (2003),  no. 1, 28--73.

\bibitem[S1]{S1} K. Saito, Extended affine root systems 1 (Coxeter transformations),
Publ. RIMS, Kyoto Univ. \textbf{21} (1985) 75-179.

\bibitem[S2]{S2} K. Saito, Extended affine root systems 2 (flat invariants),
Publ. RIMS, Kyoto Univ. \textbf{26} (1990) 15-78.

\bibitem[Se]{Se} G.B. Seligman, {\it Rational Methods in  Lie Algebras},
Lect. Notes in Pure and Applied Math. \textbf {17} Marcel Dekker, 
New York 1976.

\bibitem[Sl]{Sl} P. Slodowy, A character approach to Looijenga's invariant
theory for generalized root systems, Compositio Math. \textbf{55} 
(1985), 3-32. 

\bibitem[Y1]{Y1} Y.~Yoshii, 
Coordinate algebras of extended affine Lie algebras of type $\hbox{\rm A}_1$
J. Algebra  \textbf{234}
(2000) 128--168.

\bibitem[Y2]{Y2} Y.~Yoshii, Root-graded Lie algebras with compatible
grading, Comm. Algebra \textbf{29} (2001), 3365--3391.

\bibitem[Y3]{Y3} Y.~Yoshii, Classification of quantum tori with involution,
Canad. Math. Bull.
\textbf{45} (4), 
(2002)  711--731.

\bibitem[Y4]{Y4} Y.~Yoshii, Classification of division
{$\mathbb Z\sp n$}-graded alternative algebras, J. Algebra 
\textbf{256} (2002), 28--50.
 
\bibitem[Y5]{Y5} Y.~Yoshii, Root systems extended by an abelian group
and their Lie algebras,  J. Lie Theory  \textbf{14} (2004), 
371--394.

\bibitem[Y6]{Y6} Y.~Yoshii, Lie tori -- A simple characterization of 
extended affine Lie algebras, preprint 2003.  
 
 
\end{thebibliography}
\end{document}